\magnification=\magstep1\input amstex\nopagenumbers
\UseAMSsymbols
\baselineskip 14pt plus 2pt
\font\smc=cmcsc10
\overfullrule=0pt
\def\cal{\Cal}
\def\bfx{\text{\bf x}}
\def\bfu{\text{\bf u}}
\def\tr{\text{\pe tr}\,}
\def\bfy{\text{\bf y}}
\def\bfp{\text{\bf p}}
\def\bfz{\text{\bf z}}
\def\bfv{\text{\bf v}}
\font\f=cmsy10
\def\darrow{\buildrel\cal D\over\longrightarrow}
\def\asarrow{\buildrel a.s.\over\longrightarrow}
\def\parrow{\buildrel i.p.\over\longrightarrow}
\font\pe=cmss10
\font\bigfont=cmr10 scaled\magstep1
\def\exp{\text{\pe E}}
\def\P{\text{\pe P}}

\def\real{\Bbb R}

\def\exp{\text{\pe E}}
\def\darrow{\rightarrow_{{\f D}}}
\def\iparrow{\buildrel i.p.\over\longrightarrow}
\def\asarrow{\buildrel a.s.\over\longrightarrow}

\vskip 1in
{\bigfont
\centerline{Weak Convergence of a Collection of Random}
\smallskip
\centerline{Functions Defined by the Eigenvectors of Large}
\smallskip
\centerline{Dimensional Random Matrices}}
\smallskip
\vskip .25in
\centerline{by}
\vskip .25in
\centerline{Jack W. Silverstein}
\centerline{\sl Department of Mathematics Box 8205}
\centerline{\sl North Carolina State University}
\centerline{\sl Raleigh, NC 27605-8205, USA}
\vskip4.5in

\smallskip
\flushpar Key Words and Phrases. Weak convergence on $D[0,1]$  Haar measure, eigenvectors of random matrices, Brownian
bridge, Haar measure.

\newpage

\centerline{Abstract}\vskip.1in
For each $n$, let $U_n$ be Haar distributed on the group of $n\times n$ unitary matrices. Let $\bfx_{n,1},\ldots,\bfx_{n,m} $ denote orthogonal nonrandom unit vectors in ${\Bbb C}^n$ and let 
$\text{\bf u}_{n,k}=(u_k^1,\ldots,u_k^n)^*=U_n^*\text{\bf x}_{n,k}$, $k=1,\ldots,m$. Define the following functions on [0,1]:
$X^{k,k}_n(t)=\sqrt n\sum_{i=1}^{[nt]}(|u_k^i|^2-\tfrac1n)$, $X_n^{k,k'}(t)=\sqrt{2n}\sum_{i=1}^{[nt]}\bar u_k^iu_{k'}^i$, $k<k'$.
Then it is proven that $X_n^{k,k},\Re X_n^{k,k'}$, $\Im X_n^{k,k'}$, considered
as random processes in $D[0,1]$, converge weakly, as $n\to\infty$, to $m^2$ independent copies of Brownian bridge.  

The same result holds for the $m(m+1)/2$ processes in the real case, where $O_n$ is real orthogonal Haar distributed and  $\bfx_{n,i}\in{\Bbb R}^n$, with
$\sqrt n$ in $X^{k,k}_n$ and $\sqrt{2n}$ in  $X_n^{k,k'}$ replaced with $\sqrt{\frac n2}$ and $\sqrt{n}$, respectively.   This latter 
result will be shown to hold for the matrix of eigenvectors of $M_n=(1/s)V_nV_n^T$ where $V_n$ is $n\times s$ consisting of the entries of 
$\{v_{ij}\},\  i,j=1,2,\ldots$, i.i.d. standardized and symmetrically distributed, with each $\bfx_{n,i}=\{\pm1/\sqrt n,\ldots,\pm1/\sqrt n\}$, and $n/s\to y>0$ as $n\to\infty$.  This result extends the result in J.W. Silverstein  {\sl Ann. Probab. \bf18} 1174-1194.

These results are applied to the detection problem in sampling random vectors mostly made of noise and detecting whether the sample includes a nonrandom vector.  The matrix $B_n=\theta\bfv_n\bfv_n^*+S_n$ is studied where $S_n$ is Hermitian or symmetric and nonnegative definite with either its matrix of
eigenvectors being Haar distributed, or $S_n=M_n$, $\theta>0$ nonrandom, and $\bfv_n$ is a nonrandom unit vector.  Results are derived on the distributional behavior of the inner product of vectors orthogonal to $\bfv_n$ with the eigenvector 
associated with the largest eigenvalue of $B_n.$

\newpage

\pageno=1
\footline={\hss\tenrm\folio\hss}
\overfullrule=0pt
\tolerance =10000

 {\bf 1. Introduction} Let $\{v_{ij}\}$, $i,j=1,2,\ldots$ be i.i.d. real valued standardized 
 random variables with finite fourth moment, and for each $n$ let
${M_n=\frac1sV_nV_n^T}$, where $V_n=(v_{ij})$, $i=1,2,\ldots,n$,
$j=1,2,\ldots,s=s(n)$, and $n/s\rightarrow y>0$ as
$n\rightarrow\infty$. This paper is essentally an extension of results in [16],  where it is shown that random elements in 
$D[0,1]$, the space of r.c.l.l. function on $[0,1]$ embodied with the Skorohod metric, defined by the eigenvectors of $M_n$ converge weakly to Brownian bridge under the
assumption $v_{i\,j}$ is symmetrically distributed.
 Specifically, denote by $O_n\Lambda_n O_n^T$ the spectral 
decomposition of $M_n$, where the eigenvalues of $M_n$ are arranged
along the diagonal of $\Lambda_n$ in nondecreasing order, and the columns of the orthogonal 
matrix $O_n$, are the corresponding eigenvectors (a unique determination of $O_n$ is outlined in Section 2 of [16]).    For each $n$ let $\bfx_n\in\Bbb R^n$ be a nonrandom unit vector, and let 
$\bfy_n=(y_1,y_2,\ldots,y_n)^T=O_n^T\bfx_n$.   Define for $t\in [0,1]$
$$X_n(t)\equiv\sqrt{\tfrac n2}
\sum_{i=1}^{[nt]}(y_i^2-\tfrac1n)\quad\quad([a]\equiv\text{greatest
integer}\leq a).\leqno(1.1)$$

The main result in [16] is that when $v_{i\,j}$ is symmetrically distributed, for $\bfx_n=(\pm\frac1{\sqrt n},\pm\frac1{\sqrt n},\ldots,
\pm\frac1{\sqrt n})^T$.
$$X_n\darrow W^\circ\quad\text{ as }
n\rightarrow\infty\leqno(1.2)$$
({$D$} denoting weak convergence in $D[0,1]$) where $W^\circ$ is
Brownian bridge ([3], p. 64).  

This result is a partial answer to the question of how the matrix of eigenvectors of $M_n$ are related to the Haar measure on the group ${\cal O}_n$ 
of $n\times n$ orthogonal matrices, which occurs when $v_{11}$ is mean 0 Gaussian, That is, when $M_n$ is a matrix of Wishart type.   The question is originally raised in [13] where it is conjectured that for arbitrary centered $v_{11}$ the distribution of $O_n$ in ${\cal O}_n$ is near in some way to the Haar measure ([13][14],[15],[16], see also [12]).   This resulted in [13] to an investigation in the behavior of (1.1).  When $O_n$ is Haar distributed $\bfy$ is uniformly distributed over the unit 
sphere in $\Bbb R^n$, being the same as the normalized vector, $(\zeta_1,\ldots,\zeta_n)^T$, of i.i.d mean-zero Gaussian 
entries. 
(1.1) can then be written as
$$X_n(t)=\frac{\sqrt n}{\sqrt 2}\left(\frac{\sum_{i=1}^{[nt]}\zeta_i^2}{\sum_{i=1}^n\zeta_i^2}-\frac{[nt]}n\right)
=\frac n{\sum_{i=1}^n\zeta_i^2}\frac1{\sqrt 2}\frac1{\sqrt n}\left(\sum_{i=1}^{[nt]}(\zeta_i^2-1)-\frac{[nt]}n\sum_{i=1}^n(\zeta_i^2-1)\right).\leqno(1.3)$$
Using the fact that the fourth moment of a standard normal random variable is 3, we apply Donsker's theorem ([3], Theorem 16.1) along with standard results on weak convergence of random functions on $D[0,1]$ to arrive at  (1.2).  

In [14] and [15] it is shown that a necessary condition for (1.2) to hold for all unit vectors $\bfx_n$ is that when $\exp(v_{i\,1}^2)=1$ we must have 
$\exp(v_{1\,1}^4)=3$.   Indeed, it is shown in [15] that when $\exp(v_{11}^2)=1$ but $\exp(v_{11}^4)\neq3$, there exist sequences $\{\bfx_n\}$ of unit vectors such that $\{X_n\}$ fails to converge weakly.   This result suggests a strong relationship needs to exist between the distribution of $v_{1\,1}$ and Gaussian in order for (1.2) to hold for all sequences of unit vectors, and leaves open the possibility that this is true only when $v_{11}$ is Gaussian.

However, the result in [16] indicates some similarity of the distribution of $O_n$ to Haar measure, at least when $v_{ij}$ is symmetrically distributed and the entries of $\bfx_n$ are equally weighted.  

In this paper another property of the Haar measure on $\cal{O}_n$ is derived and is shown to be true for $v_{11}$ symmetrically distributed and on unit vectors considered in [16].    In order to provide a more complete setting, the property is stated and derived on ${\cal U}_n$, the group of $n\times n$ unitary matrices.   The corresponding statements and steps in the verification for the real case will be specified in the proof.
 
Let for $d\ge2$ an integer, and $b\ge1$, $D_d^b=\Pi_{i=1}^dD[0,b]$, and $\text{\f T}_d^{\,\,b}$ denote the smallest $\sigma$-field on $D_d^b$ in which convergence of elements in $D_d^b$ is equivalent to component-wise convergence.  We will prove the following:

{\smc Theorem 1.1} For each $n$, let $U_n$ be Haar distributed on  $\cal{U}_n$. Let $\bfx_{n,1},\ldots,\bfx_{n,m} $ denote orthogonal nonrandom unit vectors in ${\Bbb C}^n$ and let 
$\text{\bf u}_{n,k}=(u_k^1,\ldots,u_k^n)^*=U^*\text{\bf x}_{n,k}$, $k=1,\ldots,m$. Define the following functions on [0,1]:
$$X^{k,k}_n(t)=\sqrt n\sum_{i=1}^{[nt]}(|u_k^i|^2-\tfrac1n),\qquad X_n^{k,k'}(t)=\sqrt{2n}\sum_{i=1}^{[nt]}\bar u_k^iu_{k'}^i\quad k<k'\leqno(1.4)$$
(``$\bar{\,\,\,\,\,}$" denoting complex conjugate).
Then $X_n^{k,k},\Re X_n^{k,k'},\Im X_n^{k,k'}\  k<k'$, considered
as random processes in $D[0,1]$, converge weakly  in $D_{m^2}^1$ to independent copies of Brownian bridge.  

The fact that $X_n^{k,k}$ converges weakly to $W^\circ$ follows along the same lines as in (1.2) where now we use the fact that 
a vector uniformly distributed on the unit sphere in $\Bbb C^n$ can be achieved by normalizing an i.i.d. vector, 
$(z_1,\ldots,z_n)^T$, where each $z_i$ is standard complex normal (real and imaginary parts i.i.d. $N(0,1/2)$), and subsequently $\exp|z_1|^2=1$,  $\exp|z_1|^4=2$.  The reason why $\Re X_n^{k,k'},\Im X_n^{k,k'}\  k<k'$ converge weakly to $W^\circ$ will be seen in the proof.   It follows from how the proof is approached, by creating the $\text{\bf u}_{n,k}$ after applying  the Gram-Schmidt orthogonalization process on a matrix of i.i.d. standard complex Gaussians, resulting in a Haar distributed unitary matrix.  

The real case is stated in the following

{\smc Theorem 1.2} For each $n$, let $O_n$ be Haar distributed on $\cal{O}_n$.  Let $\bfx_{n,1},\ldots,\bfx_{n,m} $ denote orthogonal nonrandom unit vectors in ${\Bbb R}^n$ and let 
$\text{\bf y}_k=(y_{k,1},\ldots,y_{k,n})^T=O_n^T\text{\bf x}_{n,k}$, $k=1,\ldots,m$.  For each of these $k$ define $X_n^k$, a random element in $D[0,1]$ to be (1.1) with $y_i$ replaced 
with $y_{k,i}$  For $1\leq j<k\leq m$ define $Y_n^{jk}$, a random element of $D[0,1]$, to be
$$Y_n^{jk}(t)=\sqrt{n}\sum_{i=1}^{[nt]}y_{j,i}y_{k,i},\tag1.5$$
 Then the random functions $X_n^k,Y_n^{jk}$, $1\leq j<k\leq m$ converge weakly in $D_d^1$, $d=m(m+1)/2$, to independent Brownian Bridges. 
 
 The extension of the result in [16] is the following:
 
 {\smc Theorem 1.3} Assume $v_{11}$ is symmetrically distributed about 0, $\exp v_{11}^4<\infty$, and the $m$ orthogonal vectors $\bfx_{n,k}=(\pm1/\sqrt n,\ldots,\pm1/\sqrt n)^T$ (this of course necessitates the $n$'s to be restiricted to multiples of $2^m$).  Then, with 
 $O_n$ being the orthogonal matrix of eigenvectors of $M_n=\frac1sV_nV_n^T$, the conclusion of Theorem 1.2 holds.
 
 The motivation behind studying these quantities is to analyze the detection problem in sampling random vectors mostly made of
 noise, and determining whether the sample includes multiples of a nonrandom vector.   For example, reading off the values a bank of antennas is receiving at discrete intervals of time.   If the values consist of pure Gaussian noise, then the matrix forming the sample correlation matrix $S_n$ is modeled by a Wishart matrix, and its matrix of eigenvectors would be Haar distributed, either in $\cal O_n$ or $\cal U_n$.   Suppose at certain periods of time multiples of a nonrandom unit vector $\bfv_n$ appear, resulting in the matrix
 $$B_n=\theta\bfv_n\bfv_n^*+S_n\quad\theta>0\quad\text{nonrandom}.\tag1.6$$
 It is straightforward to verity that $\lambda_n^1$, the largest eigenvalue of $B_n$, is the unique value which solves
 $$\bfv_n^*(\lambda I-S_n)^{-1}\bfv_n=1/\theta\quad\text{for }\lambda>\lambda_{\max}(S_n)\tag1.7$$
 where $I$ is the $n\times n$ identity matrix and $\lambda_{\max}(S_n)$ is the largest eigenvalue of $S_n$.  Moreover, 
 a multiple of the corresponding eigenvector is
 $$(\lambda_n^1I-S_n)^{-1}\bfv_n.\tag1.8$$
 The goal is to understand the random behavior of this largest eigenvector for $n$ large in order to infer as much as possible the nature of $\bfv_n$.  We will place $S_n$ in a more general setting.
 
 Let, for each $n$, $S_n$ be a Hermitian nonnegative definite random matrix whose matrix of eigenvectors is Haar distributed in $\cal U_n$. Let
 $F_n$ denote the empirical distribution function of the eigenvalues of $S_n$, that is, for $x\ge0$, $F_n(x)=\frac1n(\text{number of eigenvalues of $S_n\ \leq x$})$.  Suppose with probability one $F_n$ converges in distribution to $F$, a nonrandom probability distribution function, continuous on $[0,\infty)$, where the largest eigenvalue of $S_n$ converges almost surely to $\lambda_{\max}>0$.
 
 We will prove the following:
 
  {\smc Theorem 1.4} Suppose for all $\lambda>\lambda_{\max}$, $\int(\lambda-x)^{-1}dF(x)\ \leq1/\theta$ (integral being over $[0,\lambda_{\max}]$).  Then with probability
 one $\lambda_n^1\to\lambda_{\max}$ as $n\to\infty$ and knowledge of the limiting behavior of (1.8) is beyond the scope of this paper.
 
 However, if there exists $\lambda>\lambda_{\max}$ such that $\int(\lambda-x)^{-1}dF(x)>1/\theta$, then, since $\int(\lambda-x)^{-1}dF(x)$ decreases to zero, there exists a unique $\lambda_1>\lambda_{\max}$ such that
 $$\int(\lambda_1-x)^{-1}dF(x)=1/\theta,\tag1.9$$
 and  $\lambda_n^1\asarrow\lambda_1$.   
 
For any $\bfx_{n,1},\ldots,\bfx_{n,m-1}$ unit vectors orthogonal to $\bfv_n$ 
$$\multline \sqrt{2n}\bfx_{n,k}^*(\lambda_n^1I-S_n)^{-1}\bfv_n\hfill\\ \hfill\darrow \int(\lambda_1-x)^{-1}dW_{k,r}^0(F(x))+i \int(\lambda_1-x)^{-1}dW_{k,i}^0(F(x)),\endmultline\tag 1.10$$
where $W^0_{k,r}$, $W^0_{k,i}$, $k\leq m-1$, are independent copies of Brownian bridge, and $I_A$ is the indicator function on the set $A$ .  Thus the limits are iid mean zero Gaussians, and it is straightforward to show their common variance is
$$\int
(\lambda_1-x)^{-2}dF(x)-
\left(\int(\lambda_1-x)^{-1}dF(x)\right)^2.\tag1.11$$
Moreover, the norm of the eigenvector (1.8) satisfies
$$\|(\lambda_n^1I-S_n)^{-1}\bfv_n\|\asarrow\left(\int(\lambda_1I-x)^{-2}dF(x)\right)^{1/2}.\tag1.12$$

With Theorems 1.2 and 1.3 come the analogous results in the real case, with (1.10) becoming
$$\sqrt{n}\bfx_{n,k}^*(\lambda_n^1I-S_n)^{-1}\bfv_n\darrow\int(\lambda_1-x)^{-1}dW_{k}^0(F(x)).\tag1.13$$
For the matrix $S_n=M_n$ in Theorem 1.3 the vectors $\bfv_n$ and $\bfx_{n,i}$ are all
orthonormal  vectors of the form $(\pm1/\sqrt n,\ldots,\pm1/\sqrt n)^T $.  There is a limiting $F$ in this case, described below.

These results can aid in detecting the presence of a particular signal by establishing the distributional behavior of inner products 
of the eigenvector of $B_n$ associated with the largest eigenvalue with vectors orthogonal to $\bfv_n$.  Knowledge of eigenvalue behavior of $S_n$ can aid in the detection.  For example, if $S_n=M_n$ where the $v_{ij}$ are $N(0,1)$, $F_y$ is known to be the Mar\v cenko-Pastur distribution ([10], [7], [18], [8], [19], [17]), proven in [19] under the assumption of finite second moment of $v_{11}$, where, with $a=(1-\sqrt y)^2$
$b=(1+\sqrt y)^2$, for $y\leq1$,  $F_y$ has
density

$$f_y(x)=\cases\frac{\sqrt{(x-a)((b-x)}}{2\pi yx}&a<x<b\\0&\text{otherwise,}\endcases$$ 
and, for $y>1$, $F$ has mass $1-1/y$ at 0, and
density $f_y(x)$ on $((1-\sqrt y)^2,(1+\sqrt y)^2)$.

These results have connections to the spike model ([1], [2], [11])
where a sample covariance matrix is studied with several of its population eigenvalues being altered, not enough of them to change the limiting empirical spectral distribution, but enough of a change in values to reveal individual sample eigenvalues associated with them.   For $B_n$ the size of $\theta$ in relation to the function $f(\lambda)=\int(\lambda-x)^{-1}dF(x)$ on $(\lambda_{\max},\infty)$  
determines whether a spike sample eigenvalue is revealed.  

The next sections contain proofs of these results.  Section 2 contains the proofs of Theorems 1.1 and 1.2,  Section 3 has the 
proof of Theorem 1.3, and Section 4 has the proof of Theorem 1.4

 {\bf 2. Proofs of Theorem 1.1 and 1.2.} We concentrate on the proof of Theorem 1.1 and indicate the analogous results in the real case.

We begin with understanding the relationship between $u_{n,k}$ and $u_{n,k'}$\ $k\neq k'$.  
Let $U$ be any unitary matrix having $\bfx_{n,k}\ k\leq m$ for its first
$m$ columns.  We know that the matrix $U_n^*U$ is also Haar distributed,
so we see that $\bfu_{n,k}\ k\leq m$, have the same distribution as the first
$m$ columns of a Haar distributed matrix.  The following lemma will
enable us to express their relationship in a simple way.

{\smc Lemma 2.1.} Let $Z=(z_{ij})$ be $n\times n$ consisting of i.i.d. complex
Gaussian entries ($z_{11}=z_r+iz_i$ $z_r,z_i$ independent $N(0,1/2)$).
Form the $n\times n$ unitary matrix $U$ by performing the Gram-Schmidt
process on the columns of $Z$.  Then $U$ is Haar distributed in 
$\text{\f U}_n$, the group of $n\times n$ unitary matrices. 

Proof: Let $\bfz_k,\bfu_k$ be the $k-th$ column of $Z,U$, respectively.
Then 
$$\bfu_1=f_1(\bfz_1)\equiv\left(\frac1{\|\bfz_1\|}\right)\bfz_1,$$
and recursively
$$\bfu_k=f_k(\bfz_1,\ldots,\bfz_k)\equiv\frac1{\|\bfz_k-\bfp_k\|}(\bfz_k-\bfp_k),$$
where
$$\bfp_k\equiv(\bfu_1^*\bfz_{k})\bfu_1+\cdots+(\bfu_{k-1}^*\bfz_k)\bfu_{k-1}.$$

Let $Q\in\text{\f U}_n$.  We will show for $k=1,\ldots,n$
$$Q\bfu_k=Qf_k(\bfz_1,\ldots,\bfz_k)=f_k(Q\bfz_1,\ldots,Q\bfz_k).\tag2.1$$
We use induction.  $k=1$ is obvious.  Assume it is true for
$\ell=1,2,\ldots,k-1$.  Then
$$Q\bfu_k=\frac1{\|Q\bfz_k-Q\bfp_k\|}(Q\bfz_k-Q\bfp_k),$$
and $$Q\bfp_k=((Qf_1(\bfz_1))^*Q\bfz_k)Qf_1(\bfz_1)+\cdots+((Qf_{k-1})^*Q\bfz_k)
Qf_{k-1}(\bfz_1,\ldots,\bfz_{k-1})$$
$$=((f_1(Q\bfz_1))^*Q\bfz_k)f_1(Q\bfz_1)+\cdots+((f_{k-1}(Q\bfz_1,\ldots,Q\bfz_{k-1}))^*Q\bfz_k)
f_{k-1}(Q\bfz_1,\ldots,Q\bfz_{k-1}),$$
by the inductive hypothesis.   Therefore we get (2.1).

We use now the fact that $QZ\sim Z$ to conclude 
$$QU=(Qf_1(\bfz_1),Qf_2(\bfz_1,\bfz_2),\ldots,Qf_n(\bfz_1,\ldots,\bfz_n))$$
$$=(f_1(Q\bfz_1),f_2(Q\bfz_1,Q\bfz_2),\ldots,f_n(Q\bfz_1,\ldots,Q\bfz_n))$$
$$\sim(f_1(\bfz_1),f_2(\bfz_1,\bfz_2),\ldots,f_n(\bfz_1,\ldots,\bfz_n))=U,$$
and we are done.

We will use Lemma 2.1 after we establish the framework for considering
the $m^2$ processes on a common probability space.

We assume the reader is familiar with the basic concepts of 
probability, including:
the notion of a measure space $\{\Omega,\text{\f F}\}$, where {\f F}
is a $\sigma$-field of subsets of $\Omega$, and a probabilty space
$\{\Omega,\text{\f F},\P\}$, where $\P$ is a probability measure defined 
on {\f F}. Given two measurable spaces $\{\Omega_1,\text{\f F}_1\}$,
$\{\Omega_2,\text{\f F}_2\}$, a mapping $T:\Omega_1\rightarrow\Omega_2$, is
{\sl measurable} $\text{\f F}_1/\text{\f F}_2$ 
if $T^{-1}A_2=\{\omega\in\Omega_1:T\omega\in A_2\}\in\text{\f F}_1$ for 
each $A_2\in
\text{\f F}_2$. For any collection {\f A} of subsets of a set $\Omega$,
$\sigma(\text{\f A})$ denotes the smallest $\sigma$-field containing
{\f A}.

We also assume the reader is also familiar with the material in 
 [3],[5]
on weak convergence of probability measures on metric spaces, most notably
the metric space $D=D[0,1]$ consisting of real valued functions on $[0,1]$
that are right continuous with left-hand limits, the $\sigma$-field {\f D},
defined by the Skorohod topology on $D$.  For $0\leq t_1<\cdots<t_k\leq1$, let
$\pi_{t_1\cdots t_k}$ denote the natural projection from $D$ to $\real^k$:
$$\pi_{t_1\cdots t_k}(x)=(x(t_1),\ldots,x(t_k)),$$
for any $x\in D$.  Let $\text{\f D}_f$ denote the collection, 
$\pi^{-1}_{t_1\cdots t_k}H$, for any $k$, $0\leq t_1<\cdots<t_k\leq1$, and
$H\in\text{\f R}^k$, the $\sigma$-field of Borel sets in $\real^k$, called
the class of finite-dimensional sets.   In  [5] it is shown that 
$\text{\f D}_f$ is a $\pi$-system (closed under intersections) and 
$\sigma(\text{\f D}_f)=\text{\f D}$.  Therefore (Theorem 3.3 of 
[4] $\text{\f D}_f$ is a separating class for probability measures on
$(D,\text{\f D})$:  if probability measures $\P_1$, $\P_2$ agree on 
$\text{\f D}_f$
then they are identical.   Thus, showing weak convergence of a sequence, 
$\{P_n\}$, of probability measures on $(D,\text{\f D})$ to a probability measure
$P$ (denoted by $\P_n\Rightarrow \P$) amounts to verifying $\{\P_n\}$ is tight
(that is, for any $\epsilon>0$ there exists a compact set $A_{\epsilon}\in
\text{\f D}$ such that $\P_n(A_{\epsilon})>1-\epsilon$ for all $n$), and 
$\P_n(A)\to \P(A)$ for all $A\in\text{\f D}_f$.  

We wish to extend this criterion of weak convergence to the product
space $D_d=\Pi_{i=1}^dD$ with the product topology $\text{\f T}_d$, the smallest 
$\sigma$-field in which convergence of elements in $D_d$ is equivalent
to component-wise convergence.  Since $(D,\text{\f D})$ is separable,
it follows from natural extensions to the material in M10 of 
[5], $(D_d,\text{\f T}_d)$ is separable, which implies
$$\text{\f T}_d=\sigma(\{\Pi_{i=1}^dA_i: 
\text{ each }A_i\in\text{\f D}\}).\tag2.2$$ 
Let $B=\{\Pi_{i=1}^d A_i:
\text{ each }A_i\in\text{\f D}_f\})$. It is clear that $B$ is also a 
$\pi$-system. We also have 

{\smc Lemma 2.2.}  $\sigma(B)=\text{\f T}_d.$

Proof: We have $\sigma(B)\subset\text{\f T}_d$. 
Let $T_1(x_1,\ldots,x_d)=x_1,$ and define 
$$C=\{A\in\text{\f D}:T_1^{-1}A\in\sigma(B)\}.$$
We have  obviously $D\in C$, and $A\in C$ for each 
$A\in\text{\f D}_f$, since 
$T_1^{-1}A=A\otimes\Pi_{i=1}^{d-1}D\in\sigma(B)$.  For $A\in C$, 
$T_1^{-1}A^c=(T_1^{-1}A)^c\in\sigma(B)$, which implies 
$A^c\in C$.  For $\{A_n\}\subset C$,  
$T_1^{-1}\cup A_n=\cup T_1^{-1}A_n\in\sigma(B)$, implying $\cup A_n\in C$.
Therefore, $C$ is a $\sigma$-field containing $\text{\f D}_f$, and
hence contains $\text{\f D}=\sigma(\text{\f D}_f)$. Therefore,
$C=\text{\f D}$, and we have for any $A\in\text{\f D}$ $A\otimes\Pi_{i=1}^{d-1}\in\sigma(B)$.
Similarly, we have for $2\leq j<d$ $(\Pi_{i=1}^{j-1}D)\otimes A\otimes(\Pi_{i=1}^{d-j}D)$ and $(\Pi_{i=1}^{d-1}D)\otimes A$ all contained in
$\sigma(B)$, so it also contains all $\Pi_{i=1}^d A_i$ for each 
$A_i\in\text{\f D}$.  Therefore by (2.2) we have $\text{\f T}_d\subset\sigma(B)$,
and we have our result.

We see then that from Theorem 3.3 of [4] $B$ is a
separating class for probability measures on $(D_d,\text{\f T}_d)$.

It is straightforward to verify that
$$B=\{\Pi_{i=1}^dA_i: A_i=\pi^{-1}_{t_1,\ldots,t_k}H_i,k=1,2,...\,,
0\leq t_1<\cdots<t_k\leq1,H_i\in\text{\f R}^k\}.\tag2.3$$
Suppose now we have a probability space $(\Omega,\text{\f F},\P)$ and a
mapping $X$ from $\Omega$ into $D_d$, for which each component $x_i$ is
a random element in $D$, that is, it is measurable {\f F}/{\f D}.  Then
for any $A_i\in\text{\f D}$, $i=1,\ldots,d$, we have 
$$X^{-1}(\Pi_{i=1}^d A_i)=
\bigcap_{i=1}^d\{\omega:x_i(\omega)\in A_i\}\in\text{\f F}.$$
Therefore, from (2.2) and Theorem 13.1 of [4] we have that
$X$ is measurable $\text{\f F}/\text{\f D}_d$, that is, $(x_1,\ldots,x_d)$
is a random element in $D_d$.

If $x_1,x_2,\ldots,x$ are random elements from probability space $(\Omega,\text{\f F},\P)$ to $D$ ($D_d$),
we write $x_n\Rightarrow x$ to mean the measures $x_n$ induce on $D$ ($D_d$) converge
weakly to the measure on $D$ ($D_d$) induced by $x$.  Also we say $\{x_n\}$ is tight (on $D$ or $D_d$) 
if the sequence of induced measures is tight.

We then have the following:

{Lemma 2.3.}  Suppose $\{x_n^1,\ldots,x_n^d\}$ is a sequence
of random functions, each lying in $D$, defined on a common
probability space $(\Omega,\text{\f F},\P)$ .  Then, from above, for each $n$ $\{x_n^1,\ldots,x_n^d\}$
is a random element in $D_d$.   Assume each $\{x_n^i\}$ is tight.  Moreover, assume there exists a random
element $(x^1,\ldots,x^d)$ in $D_d$ for which 
$$(x^1_n(t_1),\ldots,x^1_n(t_k),
\ldots,x^d_n(t_1),\ldots,x^d_n(t_k))
\Rightarrow
(x^1(t_1),\ldots,x^1(t_k),\ldots,x^d(t_1),\ldots,x^d(t_k))$$
(weak convergence on $\real^{dk}$) for all $k$, $t_1,\ldots,t_k$. 
Then $(x_n^1,\ldots,x_n^d)\Rightarrow(x^1,\ldots,x^d)$.

Proof.  Let $\P_n^i$, $\P^i$ denote the measures the $x_n^i$, $x^i$ induce
on $D$, and $\P_{n,d}$ the measure $\{x_n^1,\ldots,x_n^d\}$ induces on
$D_d$.  Then each $\{P^i_n\}$ is tight.  Therefore for any $\epsilon>0$
there exists compact sets $A^i_{\epsilon}\in\text{\f D}$ for which
$\P_n^i(A^i_{\epsilon})>1-\epsilon/d$.  Then ([5], M6)
we have $\Pi_{i=1}^d A^i_{\epsilon}$
compact, and 
$$\multline \P_n(\Pi_{i=1}^d A^i_{\epsilon})=\P(\{\omega:x_n^i(\omega)\in A^i_{\epsilon}, i\leq d\})=
\P(\cap\{\omega:x_n^i(\omega)\in A^i_{\epsilon}\})\\1-\P(\cup\{\omega:x_n^i\in A^i_{\epsilon}\}^c)
\ge1-\sum \P_n^i({A^i_{\epsilon}}^c)\ge1-\epsilon.\endmultline$$
Therefore $\{\P_{n,d}\}$ is tight.  Since $B$ is a separating class, and
it can be expresssed as in (3), we must have $\{x_n^1,\ldots,x_n^d\}\Rightarrow\{x^1,\ldots,x^d\}$.  

We proceed to show each of $X_n^{k,k},\Re X_n^{k,k'},\Im X_n^{k,k'}\  k<k'$ converges
weakly to independent copies of Brownian bridge.

The following lemma is needed throughout the remaining arguments.

Lemma 2.4. If random variables $X_n,Y_n$ are such that $\{Y_n\}$ is tight
and $X_n\iparrow 0$, then $X_nY_n\iparrow0$.

Proof:  For $\epsilon>0$ $M>0$ we have 
$$\P(|X_n|\,|Y_n|>\epsilon)=
\P(|X_n|\,|Y_n|>\epsilon, |Y_n|>M)+\P(|X_n|\,|Y_n|>\epsilon, |Y_n|\leq M)$$
$$\leq\P(|Y_n|>M)+\P(|X_n|>\epsilon/M).$$
Therefore $\limsup_n\P(|X_n|\,|Y_n|>\epsilon)\leq \limsup_n\P(|Y_n|>M)$ which
can be made arbitrarily small.  We get our result.

Let $Z$ and $U$ be as in Lemma 2.1.  We can assume the first $m$ columns of $U$ are the
orthonormal vectors $\bfu_{n,k}$ where in the following we suppress he dependence on $n$. 
We can also assume that $Z$ and $U$ are $n\times m$.  Define $r_{jk}=\bfu_j^*\bfz_k$ for $j<k$, $r_{11}=\|\bfz_1\|$, and for $k\geq 2$, $r_{kk}=\| \bfz_k-\bfp_k\|$.
We have then $r_{11}\bfu_1=\bfz_1$, and for $k\ge2$
$$r_{kk}\bfu_k=\bfz_k-\sum_{j=1}^{k-1}r_{jk}\bfu_j.$$
Letting $R$ denote the $m\times m$ upper triangular matrix $(r_{jk})$ we obtain the $QR$ factorization
of $Z$: $Z=UR$. Letting $A=R^{-1}$ we have $U=ZA$.   We have then for each $k$
$$ \bfu_k=a_{kk}\bfz_k+\sum_{j=1}^{k-1}a_{jk}\bfz_j.\tag2.4$$

For $j<k$, $\bfu_j$ and $\bfz_k$ are independent.  Therefore 
$$E(r_{jk})=0\quad\text{and}\quad E|r_{jk}|^2=1.\tag2.5$$
Therefore above the diagonal the entries of $R$ are tight.  By the weak law of large numbers
$$\|\bfz_k\|/\sqrt n\iparrow1.\tag2.6$$
It is straightforward to verify
$$r_{kk}^2=\|x_k\|^2-\sum_{j=1}^{k-1}|r_{jk}|^2.\tag2.7$$

Therefore we have
$$\frac{r_{kk}^2}{\|\bfz_k\|^2}=1+O(1)/n,\tag2.8$$
where here and in the following $O(1)$ denotes a tight sequence of random variables.
From (2.6) and (2.8) we get
$$r_{kk}/\sqrt n\iparrow1.\tag2.9$$
We have $a_{kk}=1/r_{kk}$  and for $j<k$ $a_{jk}=R_{kj}/\det(R)$, where  $R_{kj}$ is the $kj$
cofactor of $R$:
$$R_{kj}=(-1)^{k+j}\det(M),$$
and $M=M_{kj}$ is the $(m-1)\times(m-1)$ matrix obtained by deleting the $k^{\text{th}}$ row and $j^{\text{th}}$ column of $R$.  We have $\det(R)=\Pi_{i=1}^mr_{ii}.$   For $\det(M)$ we use the Leibniz formula
$$\det(M)=\sum_{\sigma\in\text{\f S}_{m-1}}\text{sgn}(\sigma)\Pi_{i=1}^{m-1}m_{i\sigma_i},$$
where $\text{\f S}_{m-1}$ is the set of all permutations of $\{1,\ldots,m-1\}$, the sum is over the collection of all permutations $\sigma\in\text{\f S}_{m-1}$, and $\text{sgn}(\sigma)$, the signature of $\sigma$, is $1$ if the reordering of $(1,\ldots,m-1)$ given by $\sigma$ can be brought back to $(1,\ldots,m-1)$ by successively interchanging two entries an even number of times, $-1$ if an
odd number of interchanges are needed.

We see then that $a_{jk}$ can be written as a sum of $(m-1)!$ terms.  The largest term in absolute value occurs for that $\sigma$ where all $r_{ii}$ $i\neq j,k$ are included.  The remaining entry must be $r_{jk}$.
Indeed, it will lie in row $j$ of $M$, the only row of $M$ not containing an $r_{ii}$, $i\neq j,k$, and column $k-1$ of $M$ (column $k$ of $R$) the only column of $M$ not containing an $r_{ii}$, $i\neq j,k$.  
The $\sigma$ creating this term is necessarily the top row of
$$\matrix\ldots&k-1&\ldots&k-2&\ldots\\
            \ldots&j&\ldots&k-1&\ldots\endmatrix$$
except when $k=j+1$ in which case the top row is   $1\  2\  ...\ m-1$.  Here the second row is  $1\  2\  ...\ m-1$. All other numbers in the top row are in increasing order.
When $k>j+1$ it takes $k-j-1$ pairwise interchanges to bring $k-1$ to the right of $k-2$ (no interchanges when $k=j+1$).  
Therefore $\text{sgn}(\sigma)=(-1)^{k-j-1}$, and since $j+k+k-j-1=2k-1$ we have
$$a_{jk}=-r_{jk}/(r_{jj}r_{kk})+O(1)/n^{3/2}.$$
   We have
$$r_{jk}\left(\frac1{r_{jj}r_{kk}}-\frac1n\right)=\frac{r_{jk}}n\left(\frac{n}{r_{jj}r_{kk}}-1\right),$$
so from (9)
$$a_{jk}=-r_{jk}/n+o(1)/n=O(1)/n,\tag2.10$$
where here and in the following $o(1)$ denotes a sequence of random variables converging in probability to zero.  We have
$$\multline r_{jk}=\left(\bfz_j^*\bfz_k-\sum_{i=1}^{j-1}\bar r_{ij}r_{ik}\right)/r_{jj}=\bfz_j^*\bfz_k/r_{jj}+O(1)/\sqrt n\\
=\bfz_j^*\bfz_k/\sqrt n+\frac{\bfz_j^*\bfz_k}{\sqrt n}\left(\frac{\sqrt n}{r_{jj}} -1\right)+O(1)/\sqrt n.
\endmultline$$
By the Central Limit Theorem $\bfz_j^*\bfz_k/\sqrt n$ is tight.   Therefore
$$a_{jk}=-\bfz_j^*\bfz_k/n^{3/2}+o(1)/n.\tag2.11$$

Let $\|\cdot\|$ represent the sup norm on functions.   Write $\bfz_j=(z_j^1,\ldots,z_j^n)^T$. Using (2.4) we have
$$X_n^{k,k}(t)=\sqrt n\left(\sum_{i=1}^{[nt]}|a_{kk}z_k^i+\sum_{j=1}^{k-1}a_{jk}z_j^i|^2-\frac{[nt]}n\right)$$
$$\multline=\sqrt n\left(a_{kk}^2\sum_{i=1}^{[nt]}|z_k^i|^2+\sum_{i=1}^{[nt]}|\sum_{j=1}^{k-1}a_{jk}z_j^i|^2 +a_{kk}\sum_{i=1}^{[nt]}\sum_{j=1}^{k-1}\bar a_{jk}z_k^i\bar z_j^i\right.\\ \left.+a_{kk}\sum_{i=1}^{[nt]}\sum_{j=1}^{k-1}a_{jk}\bar z_k^iz_j^i-\frac{[nt]}n\right).\endmultline$$
Using Cauchy-Schwarz, Lemma 2.4, the weak Law of Large Numbers,  and (2.10) we have
$$\|\sqrt n\sum_{i=1}^{[nt]}|\sum_{j=1}^{k-1}a_{jk}z_j^i|^2\|\leq n^{3/2}\sum_{i=1}^{k-1}|a_{jk}|^2\frac1n\sum_{j=1}^{k-1}\|\bfz_j\|^2\iparrow0\tag2.12$$

We have using (2.9) and (2.10)
$$\|\sqrt na_{kk}\sum_{i=1}^{[nt]}\sum_{j=1}^{k-1}\bar a_{jk}z_k^i\bar z_j^i\| \leq
(O(1)/\sqrt n)\sum_{j=1}^{k-1}\frac1{\sqrt n}\|\sum_{i=1}^{[nt]}z_k^i\bar z_j^i\|\iparrow0,$$
since $\|\cdot\|$ is continuous on $C[0,1]$, and the real and imaginary parts of $(\sqrt{2/n})\sum_{i=1}^{[nt]}z_k^i\bar z_j^i$, each satisfying
the assumptions of Donsker's theorem ([3], Theorem 16.1), converge weakly to Wiener measure, which lies in $C[0,1]$, so
that from Theorem 5.1 of  [3] (with $h=\|\cdot\|$) $\|\frac1{\sqrt n}\sum_{i=1}^{[nt]}z_k^i\bar z_j^i\|$ is tight, and using Lemma 2.4 we get our 
result.

From (2.6), (2.7), and (2.9)  we have
$$\|\sqrt na_{kk}^2\sum_{i=1}^{[nt]}|z_k^i|^2-\frac{\sqrt n}{\|\bfz_k\|^2}\sum_{i=1}^{[nt]}|z_k^i|^2\|=
\|\bfz_k\|^2\sqrt n|a_{kk}^2-1/\|\bfz_k\|^2|=O(1)\frac{\sqrt n}{r_{kk}^2}\iparrow0$$
Therefore
$$\|X_n^{k,k}-X_n^k\|\iparrow0,$$  
where
$$X_n^k(t)=\sqrt n\left(\frac1{\|\bfz_k\|^2}\sum_{i=1}^{[nt]}|z_k^i|^2-\frac{[nt]}n\right).$$

We have 
$$X^k_n(t)=\frac{\sqrt n}{\|\bfz_k\|^2}\left(\sum_{i=1}^{[nt]}(|z_k^i|^2-1)-
\frac{[nt]}n(\|\bfz_k\|^2-n)\right)=\frac n{\|\bfz_k\|^2}h_n(W^k_n(t)),$$
where
$$W^k_n=\frac1{\sqrt n}\sum_{i=1}^{[nt]}(|z_k^i|^2-1),$$
and $h_n:D\to D$ is defined as $h_n(X)=X(t)-([nt]/n)X(1)$.  Let 
$h(X)=X(t)-tX(1)$.
  We have for any 
$X\in D$ $\|h_n(X)-h(X)\|\leq\|X(1)|\,|t-[nt]/n|\leq |X(1)|/n\to 0$.
If $X_n\to X$ in the Skorohod topology, then there exists $\{\lambda_n$\},
each increasing continous on [0,1] with $\lambda_n(0)=0$, $\lambda_n(1)=1$,
such that $\|\lambda_n(t)-t\|\to0$ and $\|X_n(t)-X(\lambda_n(t))\|\to0$.
Therefore
$$\|h_n(X_n(t))-h(X(\lambda_n(t)))\|\leq \|h_n(X_n(t))-h_n(X(\lambda_n(t)))\|
+\|h_n(X(\lambda_n(t)))-h(X(\lambda_n(t)))\|$$
$$\leq\|X_n(t)-X(\lambda_n(t))\|+|X_n(1)-X(1)|+|X(1)|\,|([nt]/n)-t|\to0.$$
Therefore the set $E$ in Theorem 5.5 of [3] is empty, and
by (9.13), Theorem 16.1 and Theorem 5.5 of [3] we have 
$h_n(W_n^k)\darrow h(W)=W^{\circ}$, $W$ denoting Wiener measure.

We have $\|X_n^k-h_n(W^k_n)\|\leq|1-n/\|\bfz_k\|^2|\max_t|h_n(W_n^k(t))|$.
By (2.6) we have $|1-n/\|\bfz_k\|^2|\iparrow0.$ Again, from Theorem 5.1 of [3]
we have $\|h_n(W_n^k)\|\darrow\|W^{\circ}\|$.  Therefore, by Lemma 2.4 we have
$$\|X_n^k-h_n(W_n^k)\|\iparrow0.$$
Therefore, $X_n^{k,k}\darrow W^{\circ}$.

For $k<k'$
$$X_n^{k,k'}(t)=\sqrt{2n}\left(\sum_{i=1}^{[nt]}(a_{kk}\bar z_k^i+\sum_{j=1}^{k-1}\bar a_{jk}\bar z_j^i)(a_{k'k'}z_{k'}^i+\sum_{j'=1}^{k'-1}a_{j'k'}z_{j'}^i)\right)$$
$$\multline=\sqrt{2n}\left(a_{kk}a_{k'k'}\sum_{i=1}^{[nt]}\bar z_k^iz_{k'}^i+a_{kk}\sum_{j'=1}^{k'-1}a_{j'k'}\sum_{i=1}^{[nt]}\bar z_k^iz_{j'}^i\right.\\ \left.+a_{k'k'}\sum_{j=1}^{k-1}\bar a_{jk}\sum_{i=1}^{[nt]}\bar z_j^iz_{k'}^i+\sum_{i=1}^{[nt]}\bigl(\sum_{j=1}^{k-1}\bar a_{jk}\bar z_j^i\bigr)\bigl(\sum_{j'=1}^{k'-1}a_{j'k'}z_{j'}^i\bigr)\right)\endmultline$$
From Cauchy-Schwarz and (2.12) we have
$$\multline\left\|\sqrt n\sum_{i=1}^{[nt]}\bigl(\sum_{j=1}^{k-1}\bar a_{jk}\bar z_j^i\bigr)\bigl(\sum_{j'=1}^{k'-1}a_{j'k'}z_{j'}^i\bigr)\right\|\\ \leq\left\|\left(\sqrt n\sum_{i=1}^{[nt]}
\biggl|\sum_{j=1}^{k-1}\bar a_{jk}\bar z_j^i\biggr|^2\right)^{1/2}\left(\sqrt n\sum_{i=1}^{[nt]}\biggl|\sum_{j'=1}^{k'-1}a_{j'k'}z_{j'}^i\biggr|^2\right)^{1/2}\right\|\iparrow0\endmultline$$
Similar to what was done earlier we have for $j'\neq k$ and $j\neq k'$ we have both
$$\left\|\sqrt na_{kk}a_{j'k'}\sum_{i=1}^{[nt]}\bar z_k^iz_{j'}^i\right\|\quad\text{and}\quad
\left\|\sqrt na_{k'k'}\bar a_{jk}\sum_{i=1}^{[nt]}\bar z_j^iz_{k'}^i\right\|$$
converging in probability to zero.
Also
$$\left\|\sqrt na_{kk}a_{k'k'}-\frac1{\sqrt n}\sum_{i=1}^{[nt]}\bar z_k^iz_{k'}^i\right\|
=\biggl|\frac{n}{r_{kk}r_{k'k'}}-1\biggr|\left|\frac1{\sqrt n}\sum_{i=1}^{[nt]}\bar z_k^iz_{k'}^i\right\|\iparrow0.
$$
We have using (2.11)
$$\multline\left\|\frac1{\sqrt n}z_k^*z_{k'}\frac{[nt]}n+\sqrt na_{kk}a_{kk'}\sum_{i=1}^{[nt]}|z_k^i|^2\right\|\\=
\left\|\frac1{\sqrt n}z_k^*z_{k'}\frac{[nt]}n+\sqrt na_{kk}a_{kk'}[nt]+na_{kk}a_{kk'}\frac1{\sqrt n}\sum_{i=1}^{[nt]}(|z_k^i|^2-1) \right\|\endmultline$$
$$\leq\left\|\frac1{\sqrt n}z_k^*z_{k'}\frac{[nt]}n+\sqrt na_{kk}[nt](-z_k^*z_{k'}/n^{3/2}+o(1)/n)\right\|+na_{kk}|a_{kk'}|\left\|\frac1{\sqrt n}\sum_{i=1}^{[nt]}(|z_k^i|^2-1) \right\|.$$
Since the function inside the norm of the second term converges weakly to Wiener measure, the second term converges in probability to zero.  The first term is
$$\leq\left|\frac1{\sqrt n}\bfz_k^*\bfz_{k'}\right||1-\sqrt na_{kk}|+o(1)\sqrt na_{kk}\iparrow0.$$

Therefore 
$$\left\|X_n^{k,k'}-\sqrt{\frac 2n}\left(\sum_{i=1}^{[nt]}
\bar z_k^iz_{k'}^i-\frac{[nt]}nz_k^*z_{k'}\right)\right\|\iparrow0.$$

We separate out the real and imaginary parts of the process $X_n^{k,k'}$ is
approaching.  Write $z_k=z_{kr}+iz_{ki}$, $z_{k'}=z_{k'r}+iz_{k'i}$.  Then
the real and imaginary parts of $X_n^{k,k'}$ are approaching, respectively
$$\sqrt{\frac2n}\left(\sum_{j=1}^{[nt]}(z_{kr}^jz_{k'r}^j+z_{ki}^jz_{k'i}^j)-
\frac{[nt]}n\sum_{j=1}^n(z_{kr}^jz_{k'r}^j+z_{ki}^jz_{k'i}^j)\right)=h_n(W_n^{k,k',r}(t))$$
and 
$$\sqrt{\frac2n}\left(\sum_{j=1}^{[nt]}(z_{kr}^jz_{k'i}^j-z_{ki}^jz_{k'r}^j)-
\frac{[nt]}n\sum_{j=1}^n(z_{kr}^jz_{k'i}^j-z_{ki}^jz_{k'r}^j)\right)=h_n(W_n^{k,k',i}(t))$$
where 
$$W_n^{k,k',r}(t)=\sqrt{\frac2n}\sum_{j=1}^{[nt]}(z_{kr}^jz_{k'r}^j+z_{ki}^jz_{k'i}^j)
\quad\text{and}\quad W_n^{k,k',i}(t)=\sqrt{\frac2n}\sum_{j=1}^{[nt]}
(z_{kr}^jz_{k'i}^j-z_{ki}^jz_{k'r}^j).$$
It is clear now that each of $X_n^{k,k}$, $\Re X_n^{k,k'}$, $\Im X_n^{k,k'}$
converges weakly to Brownian bridge.  In order to show they
converge weakly in $D_{m^2}$ to independent copies of $W^{\circ}$, we will
show the weak convergence of the $W_n^k$, $W_n^{k,k',r}$, $W_n^{k,k',i}$ 
to $W^k$, $W^{k,k',r}$, $W^{k,k',i}$,  independent copies of Wiener measure, using (9.13),
Theorem 5.5 (on $D_{m^2}$), and Theorem 16.1 all in [3].  

Let $W_n$ denote the $m\times m$ matrix consisting of the $W_n^k$ on the diagonal, the $W_n^{k,k',r}$ on the lower diagonal, and the $W_n^{k,k',i}$ on the upper diagonal.  Let $W$ denote an $m\times m$ matrix consisting of independent copies of Wiener measure.

We have each entry of $W_n$ is tight, satisfying the first condition of Lemma 2.3.
Choose $k$, $0\leq t_1<\cdots<t_k\leq1$.  To prove
$$(W_n(t_1),\ldots,W_n(t_k))\darrow
(W(t_1),\ldots,W(t_k))\tag2.13$$
it is sufficient to show
$$\multline(W_n(t_1),W_n(t_2)-W_n(t_1),\ldots,W_n(t_k)-W_n(t_{k-1}))\\ \darrow
(W(t_1),W(t_2)-W(t_1),\ldots,W(t_k)-W(t_{k-1})).\endmultline$$
 But the $k$ matrices $W_n(t_{\ell})-W_n(t_{\ell-1})$, where $t_0\equiv0$, are independent.
 By the natural extension to Theorem 3.2 in 
[3] it is sufficient to show each of these converges in distribution.
 We use the Cram\'er-Wold device (p. 48 of [3]).  Thus we need to prove that linear combinations of the entries of $W_n(t_{\ell})-W_n(t_{\ell-1})$
 converge in distribution to the corresponding linear combinations of the entries of $W(t_{\ell})-W(t_{\ell-1})$.
Fix $A=(a_{ij})\in\Bbb R^{m\times m}$.  Let $\circ$ denote Hadamard product on $m\times m$ matrices and let $\text{\bf 1}$ denote the $m$ dimensional column vector consisting of 1's.  Let 
$$Y=\text{\bf 1}^{\text T}\bigl(A\circ\sqrt n W_n(1/n)\bigr)\text{\bf 1}.$$
We have $\exp Y=0$ and $\exp(Y^2)=\sum_{i,j}a_{ij}^2$.  Therefore, from
the central limit theorem
$$\text{\bf 1}^{\text T}\bigl(A\circ(W_n(t_{\ell})-W_n(t_{\ell-1}))\bigr)\text{\bf 1} \darrow 
N(0,(t_{\ell}-t_{\ell-1})\sum_{i,j}a_{ij}^2),$$
the same distribution as $\text{\bf 1}^{\text T}\bigl(A\circ(W(t_{\ell})-W(t_{\ell-1}))\bigr)\text{\bf 1}$.
Therefore, by Lemma 2.3, we are done.


It is clear that the analysis carries over to the real case, so that Theorem 1.2 is true.  Indeed, when $Z$ consists of
i.i.d. standard Gaussian, we use in Lemma 2.1 the fact that for any $Q\in\cal O_n$ $QX\sim X$, and for the scaling of the 
$X_n^k$ and $Y_n^{jk}$ we have now the variance of a standard Gaussian is 1, while its fourth moment is 3.  


{\bf 3. Proof of Theorem 1.3}.  We let $F_n$ denote the empirical distribution function of $M_n$ with almost sure limiting distribution function $F_y$ specified above.    We will also use the fact  [20] that, because $\exp v_{11}^4<\infty$,   $\lambda_{\max}(M_n)$, the largest eigenvalue of $M_n$ satisfies
$$\lambda_{\max}(M_n)\to(1+\sqrt y)^2\quad\text{a.s. as $n\to\infty$}.\tag3.1$$

 We begin with two lemmas.

{\smc Lemma 3.1} Let $S$ be a metric space with $X_n,X$ random elements in $S$ and $X_n\darrow X$.
Suppose for each $n$, $\ell_n$ is a random positive integer, independent of $\{X_n\}$ such that for any positive integer $j$, $\P(\ell_n\leq j)\to0$ as $n\to\infty$.   Then $X_{\ell_n}\darrow X$.

Proof: Let $A$ be an $X$-continuity set.  For any positive integer $j$ we have
$$\multline\P(X_{\ell_n}\in A|\ell_n=j)=\P(X_{\ell_n}\in A,\ell_n=j)/\P(\ell_n=j)\hfill\\ \hfill=\P(X_{j}\in A,\ell_n=j)/\P(\ell_n=j)=\P(X_j\in A).\endmultline$$
 For $\epsilon>0$ let positive integer $M_1$ be such that $|\P(X_j\in A)-\P(X\in A)|<\epsilon/2$ for all $j\ge M_1$.  Let $M\ge M_1$ be such that $\P(\ell_n\leq M_1)<\epsilon/4$ for all $n\ge M$.  Then, using 
$$\P(X_{\ell_n}\in A)=\sum_{j=1}^{\infty}\P(X_j\in A)\P(\ell_n=j)$$
we have for all $n\ge M$
$$\multline|\P(X\in A)-\P(X_{\ell_n}\in A)|\leq\sum_{j=M_1+1}^{\infty} |\P(X\in A)-\P(X_j\in A)|\P(\ell_n=j)\hfill\\ \hfill+\sum_{j=1}^{M_1} |\P(X\in A)-\P(X_j\in A)|\P(\ell_n=j)<\epsilon/2+\epsilon/2=\epsilon.\endmultline$$
Therefore since $\epsilon$ was arbitrary we have $X_{\ell_n}\darrow X.$

{\smc Lemma 3.2} Let $S'$ and $S''$ be separable metric spaces,  with $X'$, $X'_n$ random elements of $S'$, defined on probability space $\P'$, and $X''$, $X''_n$ random elements of $S''$, defined on probability space $\P''$ and let $\P=\P'\times \P"$.   Then $\{X'_n\},X'$ and $\{X''_n\},X''$ are independent on $\P$.   Suppose $X'_n\darrow X'$, $X''_n\darrow X''$ and for each $n$ there exists a positive integer-valued function $\ell_n=\ell_n(X'_n)$ for which  the $\ell_n$ satisfy the condition in Lemma 3.1.  Then $(X'_n,X''_{\ell_n})\darrow(X',X'')$ on $\P$.

Proof: From Lemma 3.1 we have $X''_{\ell_n}\darrow X''$.   Let $A'$, $A''$ be respective $X'$, $X''$ -continuity sets.   Then for each $n$
$$\multline\P(X'_n\in A',X''_{\ell_n}\in A'')=\sum_{j=1}^{\infty}\P(X'_n\in A',X''_{\ell_n}\in A'',\ell_n=j)\hfill\\ \hfill=
\sum_{j=1}^{\infty}\P(X'_n\in A',X''_j\in A'',\ell_n=j)=\sum_{j=1}^{\infty}\P(X''_j\in A'')\P(X'_n\in A',\ell_n=j)\hfill\\ \hfill=\sum_{j=1}^{\infty}\P(X''_j\in A'')\P(X'_n\in A'|\ell_n=j)\P(\ell_n=j).\endmultline$$
$$\multline\P(X'_n\in A',X''_{\ell_n}\in A'')-\P(X'\in A')\P(X''\in A'')\hfill\\ \hfill=\P(X'_n\in A',X''_{\ell_n}\in A'')-\P(X'_n\in A')\P(X''\in A'')\hfill\\ \hfill+\P(X'_n\in A')\P(X''\in A'')-\P(X'\in A')\P(X''\in A'')\hfill\\ \hfill=\sum_{j=1}^{\infty}(\P(X''_j\in A'')-\P(X''\in A''))\P(X'_n\in A'|\ell_n=j)\P(\ell_n=j)\hfill\\ \hfill+\P(X'_n\in A')\P(X''\in A'')-\P(X'\in A')\P(X''\in A'')\endmultline$$
For $\epsilon>0$ let $M_1$ be such that for all $j\ge M_1$
$$\max(|\P(X'_j\in A')-\P(X'\in A')|,|\P(X'_j\in A'')-\P(X''\in A'')|)<\epsilon/3.$$
Let $M\ge M_1$ be such that for all $n\ge M$ $\P(\ell_n\leq M_1)<\epsilon/6.$  Then for all $n\ge M$
$$\multline|\P(X'_n\in A',X''_{\ell_n}\in A'')-\P(X'\in A')\P(X''\in A'')|\hfill\\ \hfill 
\leq\sum_{j=M_1+1}^{\infty}|\P(X''_j\in A'')-\P(X''\in A'')|\P(X'_n\in A'|\ell_n=j)\P(\ell_n=j)\hfill\\ \hfill
\\ +\sum_{j=1}^{M_1}|\P(X''_j\in A''-\P(X''\in A'')|\P(X'_n\in A'|\ell_n=j)\P(\ell_n=j)+\epsilon/3<\epsilon. \endmultline$$ 
Since $\epsilon$ is arbitrary we have the result.

Recalling $Y_n^{jk}$ in (1.5), let $Y_n=Y_n^{12}$.  Much of the following are modifications to the results in [16], with $X_n$ replaced by $Y_n$, with
some being used exactly as stated in that paper.   As in [16] some of the results make assumptions more general than what 
is needed to prove Theorem 1.2, in order to be able to use them in the future. Results in [15] will also be used and modified.

We proceed to prove Theorem 2.1 of [16] with $X_n$ replaced by $Y_n$. We also assume that $X_n^i(F_n(\cdot))\darrow W_{F_y(\cdot)}^0$ on $D[0,\infty)$ for $i=1,2$.
Let $\rho$ denote the sup metric in $C[0,1]$:
$$\rho(x,y)=\sup_{t\in[0,1]}|x(t)-y(t)|\quad\text{for }x,y\in D[0,1].$$

{\smc Theorem} 3.1.  $Y_n(F_n(\cdot))$, $X_n^i(F_n(\cdot))$, $i=1,2$ all converging weakly to $W_{F_y(\cdot)}^0$, in $D[0,\infty)$, $F_n\darrow F_y$ i.p., and $\lambda_{\max}\equiv\lambda_{\max}(M_n)\to(1+\sqrt y)^2$ i.p. $\Rightarrow$ $Y_n\darrow W^0$.

Proof: The proof of Theorem 2.1 in [16] applied to $Y_n$ remains unchanged up to the middle of p. 1179.  For fixed $M_n$ let $\lambda_{(1)}<\lambda_{(2)}<\cdots<\lambda_{(t)}$ be the $t$ distinct eigenvalues of $M_n$ with 
multiplicities $m_1,m_2,\ldots,m_t$.  For fixed eigenvalue $\lambda_{(i)}$ the corresponding $m_i$ columns of $O_n$ are distributed as $O_{n,i}O_{i}$ where $O_{n,i}$ is
$n\times m_i$ containing $m_i$ orthonormal columns from the eigenspace of $\lambda_{(i)}$, and $O_{i}$ is Haar distributed in the group of $m_i\times m_i$ orthogonal matrices,
independent of $M_n$.  The coordinates of $\text{\bf y}_1$ and $\text{\bf y}_2$ corresponding to $\lambda_{(i)}$ are respectively of the form
$$(O_{n,i}O_{i})^T\text{\bf x}_{n,1}=a_{1,i}\text{\bf w}_{1,i}\quad\text{ and }(O_{n,i}O_{i})^T\text{\bf x}_{n,2}=a_{2,i}\text{\bf w}_{2,i},$$
where $a_{1,i}=\|O_{n,i}^T\bfx_{n,1}\|$, $a_{2,i}=\|O_{n,i}^T\bfx_{n,2}\|$, and $\text{\bf w}_{1,i}=(w_{1,i}^1,w_{1,i}^2,\ldots,w_{1,i}^{m_i})^T$, 
$\text{\bf w}_{2,i}=(w_{2,i}^1,w_{2,i}^2,\ldots,w_{2,i}^{m_i})^T$ are each uniformly distributed on the
unit sphere in $\Bbb R^{m_i}$.  Write
$$(O_{n,i}O_{i})^T(\text{\bf x}_{n,1}+\text{\bf x}_{n,2})=a_{1,2,i}\text{\bf w}_{1,2,i},$$
where $a_{1,2,i}=\|O_{n,i}^T(\text{\bf x}_{n,1}+\text{\bf x}_{n,2})\|$ and $\text{\bf w}_{1,2,i}=(w_{1,2,i}^1,w_{1,2,i}^2,\ldots,w_{1,2,i}^{m_i})^T$ is uniformly distributed on the
unit sphere in $\Bbb R^{m_i}$.
We have (2.4) in [16] holding for $a_i=a_{1,i}$ and $a_{2,i}$.  Also as in (2.4) in [16] we have
$$\max_{1\leq i\leq t}\sqrt n|\text{\bf x}_{n,1}^TO_{n,i}O_{n,i}^T\text{\bf x}_{n,2}|\iparrow0.\tag3.1$$
We have (2.3) in [16] for $Y_n$ becomes
$$\rho(Y_n(\cdot),Y_n(F_n(F_n^{-1}(\cdot))))=\max\Sb 1\leq i\leq t\\1\leq j\leq m_i\endSb\sqrt n\bigg|a_{1,i}a_{2,i}\sum_{\ell=1}^jw_{1,i}^{\ell}w_{2,i}^{\ell}\biggr|.\tag3.2$$
For each $i\leq t$ and $j\leq m_i$
$$\sqrt na_{1,i}a_{2,i}\sum_{\ell=1}^jw_{1,i}^{\ell}w_{2,i}^{\ell}=\frac{\sqrt n}2\left(a_{1,2,i}^2\sum_{\ell=1}^j(w_{1,2,i}^{\ell})^2
-a_{1,i}^2\sum_{\ell=1}^j(w_{1,i}^{\ell})^2-a_{2,i}^2\sum_{\ell=1}^j(w_{2,i}^{\ell})^2\right)$$
$$=\frac{\sqrt n}2\left(\biggl(a_{1,2,i}^2-2\frac{m_i}n\biggr)\sum_{\ell=1}^j(w_{1,2,i}^{\ell})^2-
\biggl(a_{1,i}^2-\frac{m_i}n\biggr)\sum_{\ell=1}^j(w_{1,i}^{\ell})^2-
\biggl(a_{2,i}^2-\frac{m_i}n\biggr)\sum_{\ell=1}^j(w_{2,i}^{\ell})^2\right)\tag a$$
$$+\frac{\sqrt n}2\left(2\frac{m_i}n\biggl(\sum_{\ell=1}^j(w_{1,2,i}^{\ell})^2-\frac{j}{m_i}\biggr)-\frac{m_i}n\biggl(\sum_{\ell=1}^j(w_{1,i}^{\ell})^2-\frac{j}{m_i}\biggr)-\frac{m_i}n\biggl(\sum_{\ell=1}^j(w_{2,i}^{\ell})^2-\frac{j}{m_i} \biggr)\right) .\tag b$$
From (3.1) above and (2.4) in [16] we see the maximum of the absolute value of (a) over all $j\leq m_i$, $1\leq i\leq t$ converges in probability to zero.  We see that the three sums in (b) are beta distributed the same as in (b) of[ 16] p. 1180.  Therefore the same arguments leading to the convergence of (2.3) of [16] to zero in probability give us the convergence of (3.2) to zero i.p.  Therefore for $y\leq1$ we have $Y_n\darrow W^0$.  

For $y>1$, the main difference is the appearance of $Y_n(t)=Y_n^{12}$ for $t<F_n(0)+1/n$.  Let $\underline{\text{\bf x}}_{n,1}=O_{n,1}^T\text{\bf x}_{n,1}$ , $\underline{\text{\bf x}}_{n,2}=O_{n,1}^T\text{\bf x}_{n,2}$, and $o_i$ denote the $i^{\text{th}}$ column of $O_{1}$. Notice that $a_{i,1}=\|\underline{\text{\bf x}}_{n,i}\|$, $i=1,2$.  We have 
$$X_n^i(F_n(0))=\sqrt{\frac n2}(a_{i,i}^2-F_n(0))\darrow W_{F_y(0)}\quad \text{ as }n\to\infty$$
$i=1,2$. therefore, from Lemma 2.4
$$a_{i,1}^2\iparrow  F_y(0)=1-(1/y), \quad i=1,2.\tag3.3$$ Write
$$\underline{\text{\bf x}}_{n,1}=\frac{\underline{\text{\bf x}}_{n,1}^T\underline{\text{\bf x}}_{n,2}}{a_{2,1}^2}\underline{\text{\bf x}}_{n,2}+
\underline{\text{\bf z}}.$$  
We have $\underline{\text{\bf z}}^T\underline{\text{\bf x}}_{n,2}=0$ and 
$$\|\underline{\text{\bf z}}\|=\frac{\sqrt{a_{1,1}^2a_{2,1}^2-(\underline{\text{\bf x}}_{n,1}^T\underline{\text{\bf x}}_{n,2})^2}}{a_{2,1}}.$$  
Notice that $\sqrt{n}\underline{\text{\bf x}}_{n,1}^T\underline{\text{\bf x}}_{n,2}=Y_n(F_n(0))$.  Therefore from Lemma 2.4
$$\underline{\text{\bf x}}_{n,1}^T\underline{\text{\bf x}}_{n,2}\iparrow0.\tag3.4$$
 For $t<F_n(0)+1/n$  
$$\multline Y_n(t)=\sqrt{n}\sum_{i=1}^{[nt]}\underline{\text{\bf x}}_{n,1}^To_io_i^T\underline{\text{\bf x}}_{n,2}=\sqrt{\frac2{F_n(0)}}\frac{Y_n(F_n(0))}{\sqrt n}A_n(t)+\frac{\sqrt{a_{1,1}^2a_{2,1}^2-(\underline{\text{\bf x}}_{n,1}^T\underline{\text{\bf x}}_{n,2})^2}}{\sqrt{F_n(0)}}B_n(t)\hfill\\ \hfil+Y_n(F_n(0))\frac{[nt]}{nF_n(0)}\endmultline$$ 
where
$$A_n(t)=\sqrt{\frac{nF_n(0)}2}\biggl(\sum_{i=1}^{[nt]}\frac{\underline{\text{\bf x}}_{n,2}^T} 
{a_{2,1}}o_io_i^T\frac{\underline{\text{\bf x}}_{n,2}} 
{a_{2,1}}-\frac{[nt]}{nF_n(0))}\biggr)$$ 
and
$$B_n(t)=
\sqrt{nF_n(0)}\sum_{i=1}^{[nt]}\frac{\underline{\text{\bf z}}^T}{\|\underline{\text{\bf z}}\|}o_io_i^T\frac{\underline{\text{\bf x}}_{n,2}} 
{a_{2,1}}$$
Since $O_{1}$ is Haar distributed and independent of $M_n$, we see that $A_n$ and $B_n$ have the same distribution if $\underline{\text{\bf x}}_{n,2}/ 
{a_{2,1}}$ and $\underline{\text{\bf z}}/{\|\underline{\text{\bf z}}\|}$ were nonrandom orthonormal vectors.  $H_n(t)$ in [16] now becomes
$$H_n(t)=\sqrt{\frac2{F_n(0)}}\frac{Y_n(F_n(0))}{\sqrt n}A_n(F_n(0)\varphi_n(t))+\frac{\sqrt{a_{1,1}^2a_{2,1}^2-(\underline{\text{\bf x}}_{n,1}^T\underline{\text{\bf x}}_{n,2})^2}}{\sqrt{F_n(0)}}B_n(F_n(0)\varphi_n(t))$$
$$+Y_n(F_n(0))\left(\frac{[nF_n(0)\varphi_n(t)]}{nF_n(0)}-1\right)+Y_n(F_n(F_n^{-1}(t)),$$
where $\varphi_n(t)=\min(t/F_n(0),1)$ for $t\in [0,1]$.  Denote the sum of the last two terms by (a).
Notice that for $s\in [0,1]$, from Theorem 1.2, both $A_n(F_n(0)s)$ and $B_n(F_n(0)s)$ converge weakly to independent Brownian bridges.  We apply Lemma 3.2 where $X'_n=((a),F_n(0))$, $\ell_n=nF_n(0)$, and $X''_{\ell_n}=(A_n(F_n(0)s),B_n(F_n(0)s)$.  Since, from (3.3) and (3.4) the coefficient of $A_n$  converges i.p. to zero and the coefficient of $B_n$ converges i.p. to $\sqrt{F_y(0)}=\sqrt{1-(1/y)}$  we have $H_n$ converging weakly to $H$ appearing in [16] (notice the misprint on line 8, p. 1183 of [16].  The zero to the right of the arrow should be $\varphi(t)$).  The final argument is exactly the same as in [16].  This completes the proof of the theorem

 The next step is to extend Theorem 3.1 of [16] to random elements in $(D_d^b,\text{\f T}_d^b)$.  We denote the modulus of continuity of $x\in D[0,b]$ by $w(x,\cdot)$:
$$w(x,\delta)=\sup_{|s-t|<\delta}|x(s)-x(t)|, \quad \delta\in (0,b].$$

{\smc Theorem} 3.2. Let $\{(x_n^1,\ldots,x_n^d)\}$ be a sequence of random elements of $D_d^b$, defined on a common probability space, each $\{x_n^i\}$ satisfy the assumptions of Theorem 15.5 of [3]: $\{x_n^i(0)\}$ is tight and for every positive $\epsilon$ and $\eta$, there exists a $\delta\in(0,b)$ and an integer $n_0$, such that, for all $n>n_0$, $\P(w(x_n^i,\delta)\ge\epsilon)\leq\eta$.  If there exists a random element $(x^1,\ldots,x^d)$ with $\P(x^i\in C[0,b])=1$ for each $i$, and such that
$$\left\{\left(\int_0^bt^rx_n^1dt,\ldots,\int_0^bt^rx_n^ddt\right)\right\}_{r=0}^{\infty}\darrow
\left\{\left(\int_0^bt^rx^1dt,\ldots,\int_0^bt^rx^ddt\right)\right\}_{r=0}^{\infty}\quad \text{as }
n\to\infty\tag3.5$$
($D$ denoting weak convergence on $\Bbb R^{\infty}$), then $(x_n^1,\ldots,x_n^d)\Rightarrow(x^1,\ldots,x^d)$.

Proof. From Theorems 5.1 and 15.5 of [3] and Lemma 2.3  weak convergence will follow from showing the distribution of $$(x^1(t_1),\ldots,x^1(t_k),\ldots,x^d(t_1),\ldots,x^d(t_k))$$ for all $k$, $t_1,\ldots,t_k\in[0,1]$ is uniquely determined by the distribution of
$$\left\{\left(\int_0^1t^rx^1dt,\ldots,\int_0^1t^rx^ddt\right)\right\}_{r=0}^{\infty}.\tag3.6$$
This is achieved by showing the distribution of
$$\sum_{i=1}^d\sum_{j=1}^ka_{ij}x^i(t_j)$$
is uniquely determined by the distribution of (3.6).  By a simple extension of the proof of Theorem 3.1 in [16] this can be done.

Next we prove the analog of Theorem 4.2 in [16].
Write 
$$Y_n(F_n(x))=\sqrt n\text{\bf x}_{n,1}^TP^{M_n}([0,x])\text{\bf x}_{n,2},$$
$P^{M_n}(A)$ being the projection matrix on the subspace of $\Bbb R^n$ spanned by the eigenvectors of $M_n$ having eigenvalues in $A$, a measurable subset of $\Bbb R^+$. Assuming $v_{11}$ is symmetric, we have the following results from [16]:

{\smc Fact 3} in [16]: $P^{M_n}(A)\sim OP^{M_n}(A)O^t$ for any permutaion matrix $O$.

{\smc Lemma 4.1} in [16]: If one of the indices $i_1,j_1,\ldots,i_4,j_4$ appears an odd number of times, then for Borel sets $A_1,\ldots,A_4\in\Bbb R^+$
$$\exp\bigl(P^{M_n}_{i_1j_1}(A_1)P^{M_n}_{i_2j_2}(A_2)P^{M_n}_{i_3j_3}(A_3)P^{M_n}_{i_4j_4}(A_4)\bigr)=0.$$

Assume also that
each $\text{\bf x}_{n,j}=(\pm 1/\sqrt n,\ldots\pm 1/\sqrt n)^T$ and are orthogonal.  Then necessarily $n$ is even, say $n=2p$, and exactly $p$ entries of $\text{\bf x}_{n,2}$ are of opposite sign with the corresponding entries of $\text{\bf x}_{n,1}$.  Moreover, Fact 3 in [16] is true for  $O$ diagonal with $\pm 1$'s on its diagonal, using exactly the same argument.  If $O$ is diagonal of this type with signs matching those of $\text{\bf x}_{n,1}$ coordinatewise, then 
$$Y_n(F_n(x))=\sqrt n(O\text{\bf x}_{n,1})^TOP^{M_n}([0,x])O^TO\text{\bf x}_{n,2}\sim
\sqrt n(O\text{\bf x}_{n,1})^TP^{M_n}([0,x])O\text{\bf x}_{n,2}.\tag3.7$$
Therefore we can assume the sign of all the entries of $\text{\bf x}_{n,1}$ are positive.  Let now $O$ be a permutation matrix which moves all the positive entries of the new $\text{\bf x}_{n,2}$ to the first $p$ positions.  Then using (3.7) again we conclude that we can assume that all the entries of $\text{\bf x}_{n,1}$ and the first $p$ entries of $\text{\bf x}_{n,2}$ are positive, and that the remaining entries of $\text{\bf x}_{n,2}$ are negative.

{\smc Theorem 3.3.}  Assume $v_{11}$ is symmetrically distributed about 0, $\bfx_{n,j}=(\pm1/\sqrt n,\ldots,\pm1/\sqrt n)^T$, $j=1,2$, and are
orthogonal.  Then  
$$\exp\bigl(Y_n(F_n(0))\bigr)^4\leq\exp(27P_{11}^{M_n}(\{0\}))^2\tag3.8$$
and for $0\leq x_1\leq x_2$
$$\exp\bigl(Y_n(F_n(x_2))-Y_n(F_n(x_1))\bigr)^4\leq\exp(27P_{11}^{M_n}((x_1,x_2]))^2\tag3.9$$

Proof: With $A=\{0\}$ or $(x_1,x_2]$ (corresponding to (3.8), (3.9) respectively, we have
$$\multline\exp\bigl(Y_n(F_n(0))\bigr)^4=\frac1{n^2}\exp\left(\sum_{i\leq n;j\leq p}P_{ij}^{M_n}(A)-\sum_{p+1\leq i,j\leq n}P_{ij}^{M_n}(A)\right)^4\hfill\\ \hfill=\frac1{n^2}\exp\left(\sum_{i\leq p}P^{M_n}_{ii}(A)-\sum_{p+1\leq i\leq n}P^{M_n}_{ii}(A)+2\sum_{i<j\leq p}P^{M_n}_{ij}(A)-2\sum_{p+1\leq i<j\leq n}P^{M_n}_{ij}(A)\right)^4\endmultline\tag3.10$$
$\leq$ (using for nonnegative $a$, $b$, $c$ $(a+b+c)^4\leq27(a^4+b^4+c^4)$)
$$\frac{27}{n^2}\exp\left(\sum_{i\leq p}P^{M_n}_{ii}(A)-P^{M_n}_{i+p\ i+p}(A)\right)^4\tag a$$
$$+$$
$$\frac{54}{n^2}\exp\left(\sum_{\buildrel {i,j\leq p}\over {i\neq j}}P_{ij}^{M_n}(A)\right)^4,\tag b$$
where in (b) we used Fact 3 of [16], which says that $P^{M_n}$ is distributed the same as $OP^{M_n}O^T$ for permutation matrices $O$, on the $P_{ij}^{M_n}$'s with $i\neq j$ and both larger than $p$.  Suppressing the dependence on $M_n$ and $A$, we have from Fact 3 and Lemma 4.1 in [16] 
$$\multline (\text{b})=\frac{216p(p-1)}{n^2}\bigl(12(p-2)\exp(P_{12}^2P_{13}^2)+3(p-2)(p-3)\exp(P_{12}^2P^2_{34})\hfill\\ \hfill+12(p-2)(p-3)\exp(P_{12}P_{23}P_{34}P_{14})+2\exp(P_{12}^4)\bigr).\endmultline$$
Bounds involving $\exp(P_{12}P_{23}P_{34}P_{14})$ and $\exp(P_{12}^2P^2_{34})$ were derived in [16], from which we get
$$(n-2)(n-3)\exp(P_{12}P_{23}P_{34}P_{14})\leq\exp(P_{11}P_{22})$$
and
$$(n-2)(n-3)\exp(P_{12}^2P^2_{34})\leq\exp(P_{11}P_{22}).$$
A bound on $\exp(P^2_{12}P^2_{13})$ is also needed. Starting from the fact that $P^2=P$, we take the expected value of both sides of
$$P^2_{12}\biggl(\sum_{j\ge3}P^2_{1j} + P^2_{11} + P^2_{12}\biggr)=P^2_{12}P_{11}$$
and use Fact 3 in [16] to get 
$$(n-2)\exp(P^2_{12}P^2_{13})\leq\exp(P^2_{12}P_{11}).$$
Therefore for $p\ge2$
$$(\text{b})\leq\frac{216p(p-1)}{n^2}\left(12\frac{p-2}{n-2}\exp(P_{12}^2P_{11})+15\frac{(p-2)(p-3)}{(n-2)(n-3)}\exp(P_{11}P_{22})
+2\exp(P_{12}^4)\right).$$
Thus, using Fact 3 in [16] and the facts that $P_{11}\in[0,1]$, $P_{12}^2\leq P_{11}P_{22}$ since $P$ is nonnegative definite, and $ab\leq\frac12(a^2+b^2)$, we get
$$(\text{b})< 648\exp(P^2_{11}).$$


In (a) we expand the fourth power of the sum.  Using Fact 3 in [16] we see that any term involving  an odd number of $P_{ii}-P_{i+p\ i+p}$ is zero.  Therefore
$$\multline (a)=\frac{27p}{n^2}(\exp(P_{11}-P_{22})^4+3(p-1)\exp(P_{11}-P_{22})^2(P_{33}-P_{44})^2
\leq27\exp(P_{11}-P_{22})^2\hfill\\ \hfill
\leq54\exp P_{11}^2.\endmultline$$

Therefore, the expression in (3.10) is bounded by $ \exp(27P_{11})^2$, and the proof is complete.

Notice that for unit $\text{\bf x}_n\in \Bbb R^n$ $\text{\bf x}_n^TP_{11}^{M_n}(\cdot)\text{\bf x}_n$ is a (random) probability measure with mass at the eigenvalues of $M_n$.  In [15] it is proven that 
$$\{\sqrt{n/2}(\text{\bf x}_n^TM_n^r\text{\bf x}_n-(1/n)\tr(M_n^r)))\}_{r=1}^{\infty}\darrow\left\{\int_{(1-\sqrt y)^2}^{(1+\sqrt y)^2}x^rdW^0_{F_y(x)}\right\}_{r=1}^{\infty}\quad as \quad n\to\infty\tag3.11$$
(${\f D}$ denoting weak convergence on $\Bbb R^{\infty}$) for every sequence $\{\text{\bf x}_n\}$, $\text{\bf x}_n\in\Bbb R^n$, $\|\text{\bf x}_n\|=1$ if and only if $\exp v_{11}=0$, $\exp v_{11}^2=1$, and $\exp v_{11}^4=3$.  It is proven by showing the mixed moments of the left side of (3.11) depends on the first, second and fourth moment of $v_{11}$ after two sets of truncations and centralizations.  After the final truncation and centralization the mixed moments are shown to be bounded regardless of the value of the fourth moment as long as it is finite.  Thus after removing the $\sqrt n$ on the left side of (3.11) we find that the difference of the moments of the distribution $\text{\bf x}_n^TP^{M_n}(\cdot)\text{\bf x}_n$ and that of $F_n$, the empirical distribution of the eigenvalues of $M_n$, approach each other i.p. as $n\to\infty$.  Since it is known that $F_n\darrow F_y$ a.s. from the method of moments we conclude that $$\text{\bf x}_n^TP^{M_n}(\cdot)\text{\bf x}_n\darrow F_y\quad\text{ i.p.}$$

With $\bfx_n=(1,0,\dots,0)^T$ we conclude that
$$P_{11}^{M_n}(\cdot)\darrow F_y\quad\text{i.p.}\tag3.12$$

The next results extends (3.11) to several different $\text{\bf x}_n$'s simultaneously.  

{\smc Theorem} 3.4.   Assume $\exp v_{11}=0$ and $\exp v_{11}^2=1$.  Fix $d$ a positive integer.  Let for every $n$ $\bfx^1_n,\ldots,\bfx^d_n$, $\bfx^j_n=(x_1^j,\ldots,x_n^j)^T$, be $d$ unit vectors in $\Bbb R^n$.
Then the limiting distributional behavior of
$$\{\sqrt{n/2}({\bfx_n^1}^TM_n^r\bfx_n^1-(1/n)\tr(M_n^r)),\ldots,\sqrt{n/2}({\bfx_n^d}^TM_n^r{\bfx_n^d}-(1/n)\tr(M_n^r))\}_{r=1}^{\infty}\tag3.13$$
is the same as that when $v_{11}$ is $N(0,1)$ if either:
\item{a)} $\exp v_{11}^4=3$ or
\item{b)} for each $j\leq d$
$$\sum_{i=1}^{n}({x_i^j})^4\to0,\quad\text{ as }n\to\infty.$$

Proof of a). By [15], through a series of truncations and centralizations, it is sufficient to assume that $v_{ij}=v_{ij,n}$ iid with $|v_{11}|\leq2n^{1/4}$, $\exp v_{11}=0$, $\exp v_{11}^2\to1$,  $\exp v_{11}^4\to3$ as $n\to\infty$, and $(1/n)\tr M_n^r$ can be replaced by $\exp{\bfx_n^i}^T\!\!M_n^r\bfx_n^i$.  We will use the method of moments.  We will show for positive integers $m_1,\ldots,m_d$, $r_j^i$, $i\leq d$, $j\leq m_i$, with $m=\sum_{i=1}^dm_i$, the limiting behavior of
$$\multline n^{m/2}\exp[({\bfx_n^1}^T\!\!M_n^{r_1^1}\bfx_n^1-\exp{\bfx_n^1}^T\!\! M_n^{r_1^1}\bfx_n^1)\cdots({\bfx_n^1}^T\!\!M_n^{r_{m_1}^1}\bfx_n^1-\exp{\bfx_n^1}^T\!\! M_n^{r_{m_1}^1}\bfx_n^1)\hfill\\ \hfill\cdots({\bfx_n^d}^T\!\!M^{r_1^d}\bfx_n^d-\exp {\bfx_n^d}^T\!\!M^{r_1^d}\bfx_n^d)\cdots({\bfx_n^d}^T\!\!M^{r_{m_d}^d}\bfx_n^d-\exp{\bfx_n^d}^T\!\!
M^{r_{m_d}^d}\bfx_n^d)]\endmultline\tag3.14$$
depends only on $\exp v_{11}^2$ and $\exp v_{11}^4$ and therefore is the same when the original $v_{ij}$'s are $N(0,1)$.

 Let $r=\sum_{i=1}^d\sum_{j=1}^{m_i}r_j^i$.  We have 
$$(s^r/n^{m/2})\times(3.14)\qquad\qquad\qquad\qquad\qquad\qquad\qquad
\qquad\qquad\qquad\qquad\qquad\qquad\qquad\qquad\quad$$
$$\multline\!\!\!=\!\!\!\!\!\!\!\!\!\!\!\!\!\!\sum\Sb i^{1 1},j^{1 1},i_2^{1 1},\ldots,i_{r_1^1}^{1 1},
k_1^{1 1},\ldots,k_{r_1^1}^{1 1}\\ \vdots\\ i^{1 m_1},j^{1 m_1},i_2^{1 m_1},\ldots,i_{r_{m_1}^1}^{1  m_1},
k_1^{1 m_1},\ldots,k_{r_{m_1}^1}^{1 m_1}\\ \vdots\\i^{d 1},j^{d 1},i_2^{d 1},\ldots,i_{r_1^d}^{d 1},k_1^{d 1},\ldots,k_{r_1^d}^{ 1}\\ \vdots\\ i^{d m_d},j^{d m_d},i_2^{d m_d},\ldots,i_{r_{m_d}^d}^{d m_d},k_1^{d m_d},\ldots,k_{r_{m_d}^d}^{d\ m_d}\endSb\!\!\!\!\!\!\!\!\!\!\!\!\!\! x_{i^{11}}^1x_{j^{11}}^1\cdots x_{i^{1m_1}}^1x_{j^{1m_1}}^1\cdots x_{i^{d1}}^dx_{j^{d1}}^d\cdots x_{i^{dm_d}}^dx_{j^{dm_d}}^d\hfill\\ 
\hfill
\exp\bigl[\Pi_{\ell=1}^d\Pi_{\ell'=1}^{m_{\ell}}(v_{i^{\ell}k^{\ell\ell'}_1}v_{i_2^{\ell\ell'}k^{\ell\ell'}_1}
\cdots v_{j^{\ell\ell'}
k_{r^{\ell}_{\ell'}}^{\ell\ell'}}-\exp(v_{i^{\ell\ell'}k^{\ell\ell'}_1}v_{i_2^{\ell\ell'}k^{\ell\ell'}_1}
\cdots v_{j^{\ell\ell'}
k_{r^{\ell}_{\ell'}}^{\ell\ell'}}))\bigr].\endmultline\tag3.15$$

Now the only difference between (3.14) here and (3.15) of  [15] is that (3.15) in [15] involves only one unit vector whereas (3.14) here involves $d$ unit vectors.   The value $m=\sum m_i$ here, which is the total number of moments considered in (3.14), can be identified with the $m$ in [15], the number of moments considered in (3.15) of [15]. The expected value in (3.15) here is essentially the same as the expected value in (3.16) in  [15].  The dependence of the unit vector $\text{\bf x}_n$ in the argument presented in [15] is that the absolute value of the sum of its entries is bounded by $n^{1/2}$, its entries are bounded by 1 in absolute value, and its length is bounded.   The argument here is identical to the one in [15] using the additional fact that $|\sum_{i=1}^nx_i^jx_i^k|\leq1$ for $j,k\in\{1,\ldots,d\}$.  We have then a).  

Proof of b).  The proof follows exactly the same as in the proof of Theorem 4.1 in [16] using the additional fact that for $j_1,\ldots,j_4\in\{1,\ldots,d\}$
$$\sum_{i=1}^nx^{j_1}_ix^{j_2}_ix^{j_3}_ix^{j_4}_i\leq\max_{k\leq4}\sum_{i=1}^n{x_i^{j_k}}^4.$$
This completes the proof of Theorem 3.4

Notice that
$$\sqrt{n/2}({\bfx_{n,k}}^TM_n^r\bfx_{n,k}-(1/n)\tr M_n^r)=\int_0^{\infty}x^rdX_n^k(F_n(x))=
-\int_0^{\infty}rx^{r-1}X_n^k(F_n(x))dx\tag3.16$$
for $k\leq m$, and for $j<k$
$$\sqrt n{\bfx_{n,j}}^TM_n^r\bfx_{n,k}=-\int_0^{\infty}rx^{r-1}Y_n^{jk}(F_n(x))dx.\tag3.17$$

When $v_{11}$ is $N(0,1)$ we have from Theorem 1.2 the conclusion of Theorem 1.3.  Therefore, from
Theorem 5.1 of [3] the quantities in (3.16) and (3.17) converge weakly, together with the
quantites
$$\sqrt{n/2}\left(\frac{(\bfx_{n,j}+\bfx_{n,k})^T}{\sqrt 2}M_n^r\frac{(\bfx_{n,j}+\bfx_{n,k})}
{\sqrt 2}-(1/n)\tr M_n^r\right).\tag3.18$$ 
since
$$\multline (3.18)=\frac12\sqrt{n/2}(\bfx_{n,j}^TM_n^r\bfx_{n,j}-(1/n)\tr M_n^r)+\frac12\sqrt{n/2}(\bfx_{n,k}^TM_n^r\bfx_{n,k}-(1/n)\tr M_n^r)\hfill\\ \hfill+\sqrt{n/2}\bfx_{n,j}^TM_n^r\bfx_{n,j}.\endmultline$$

Therefore, when the $m(m+1)/2$ vectors $\bfx_{n,k}$ and $\frac{(\bfx_{n,j}+\bfx_{n,k})}
{\sqrt 2}$ are considered in Theorem 3.4 and either a) or b) hold then the quantities in (3.16) and (3.18) converge weakly to random variables having the same distribution as when $v_{11}$ is $N(0,1)$.    Since the quantity in (3.17) can be written as a linear combination of quantities in (3.16) and (3.18) we conclude that when a) or b) hold the quantities
$$\int_0^{\infty}x^rX_n^k(F_n(x))dx\quad k\leq m,\quad\int_0^{\infty}x^rY_n^{jk}(F_n(x))dx\quad j<k$$
converge weakly to random variables,  the same distribution as when $v_{11}$ is N(0,1).  Using (3.1) we have, when $b>(1+\sqrt y)^2$ 
$$\int_0^{b}x^rX_n^k(F_n(x))dx\quad k\leq m,\quad\int_0^{b}x^rY_n^{jk}(F_n(x))dx\quad j<k$$
converging weakly to variables with the same distribution as when $v_{11}$ is $N(0,1)$.  Therefore, we have (3.5) of Theorem 3.2.   Under the assumptions of Theorem 3.3 we have (3.8), (3.9) , and (3.12), which can be used as in the last paragraph of 
[16] to show that the $Y_n^{jk}(F_n(\cdot))$ also satisfy the assumptions of Theorem 15.5 of [3].  Therefore, under the assumptions of Theorem 1.3, from Theorem 1.2 and Theorem 3.2, for each $b>(1+\sqrt y)^2$ we have the $X_n^k(F_n(\cdot))$, $Y_n^{jk}(F_n(\cdot))$, $j<k$ all converging weakly in $D_d^b$ to independent copies of Brownian bridge, composed with $F_y$, and hence the convergence is also on $D[0,\infty)$ for each  of the processes .   From Theorem 2.1 in [16] and Theorem 3.1 in 
this paper,  we have the $X_n^k(\cdot)$, $Y_n^{jk}(\cdot)$ each converging weakly to Brownian bridge.   The proof of Theorem 1.3 will follow once it is shown there is joint convergence to independent copies.  

Notice that each of the limits $X_n^k(\cdot)$, $Y_n^{jk}(\cdot)$  reside in $C[0,1]$ and the limits $X_n^k(F_n(\cdot))$, $Y_n^{jk}(F_n(\cdot))$ in $C[0,\infty)$, where the topology in the latter is obtained from uniform convergence on $[0,b]$ for every $b>0$.  In fact the latter limits reside in the closed set
$$C'\equiv \{x\in C[0,\infty):x(t)=x_0\text{ for }t\in[0,(1-\sqrt y)^2] \text{ and for some } x_0,0\text{ for }t\in[(1+\sqrt y)^2,\infty)\}.$$
Consider first $y\leq1$.  Then we can assume that there is one $x_0$ in $C'$, namely 0. Let $\text{\f C}^0$ denote the class of Borel sets in $C[0,1]$ and $\text{\f C}'$ the class of Borel sets in $C'$. Define $F_y^{-1}$ to be $(1-\sqrt y)^2$ for $t=0$, $(1+\sqrt y)^2$ for $t=1$, and $F_y^{-1}(t)$ for $t\in(0,1)$.  It is straightforward to verify that the map $X(\cdot)\rightarrow X(F_y^{-1}(\cdot))$ from $C'$ to $C[0,1]$ is continuous and is the inverse of $X(\cdot)\rightarrow X(F_y(\cdot))$ from $C[0,1]$ to $C'$.  Let $\{W^{0k}_{F_y(\cdot)},W^{0jk}_{F_y(\cdot)},j<k\leq m\}$ denote the weak limit of $\{X_n^k(F_n(\cdot))$, $Y_n^{jk}(F_n(\cdot)), j<k\leq m\}$, where the entries of 
$\{W^{0k}_{(\cdot)},W^{0jk}_{(\cdot)},j<k\leq m\}=\{W^{0k}_{F_y(F_y^{-1}(\cdot))},W^{0jk}_{F_y(F_y^{-1}(\cdot))},j<k\leq m\}$ are independent copies of Brownian bridge.   Let for $A\in \text{\f C}'$ $F_y^{-1}(A)=\{X\in D[0,1]: X(F_y(\cdot))\in A\}$  be the inverse image of $A$ under $F_y^{-1}$.  Then $F_y^{-1}(A)\in\text{\f C}^0$. Suppose $A_k,A_{jk}\in\text{\f C}'$ for $j<k\leq m$.  Then
$$\multline \P(W^{0k}_{F_y(\cdot)}\in A_k,W^{0jk}_{F_y(\cdot)}\in A_{jk},j<k\leq m)\hfill\\ \hfill=\P(W^{0k}_{(\cdot)}\in F_y^{-1}(A_k),W^{0jk}_{(\cdot)}\in F_y^{-1}(A_{jk}),j<k\leq m)\hfill\\ \hfill=\Pi_kP(W^{0k}_{(\cdot)}\in F_y^{-1}(A_k))\times
\Pi_{j<k}\P(W^{0jk}_{(\cdot)}\in F_y^{-1}(A_{jk}))\hfill\\ \hfill=\Pi_k\P(W^{0k}_{F_y(\cdot)}\in A_k)\times\Pi_{j<k}\P(W^{0jk}_{F_y(\cdot)}\in A_{jk}).\endmultline\tag3.19$$
Therefore the entries of $\{W^{0k}_{F_y(\cdot)},W^{0jk}_{F_y(\cdot)},j<k\leq m\}$ are independent.

  Using the same argument used in Lemma 2.3, the sequence $\{X_n^k,Y_n^{jk},j<k\}_{n=1}^{\infty}$ is tight.
Suppose on some subsequence $\{X_n^k,Y_n^{jk},j<k\leq m\}$ converges weakly to the random element 
$\{W^{0k},W^{0jk},j<k\leq m\}$ in $D_d^1$.  Then each entry is Brownian bridge and the entries of
$\{W^{0k}_{F_y(\cdot)},W^{0jk}_{F_y(\cdot)},j<k\leq m\}$ are independent.   We invoke Theorem 8.3.7 of [6]:
 Let $X$ and $Y$ be Polish spaces (separable and can be metrized with a complete metric),  let $A$ be a Borel subset of $X$, and let $f:A\rightarrow Y$ be Borel measurable and injective (1-to-1). Then $f(A)$ is a Borel subset of $Y$.

Therefore, with $F_y(A)$ denoting the image of $A$ under $F_y$, for sets $A_k,A_{jk}\in \text{\f C}^0$ we have  $F_y(A_k),F_y(A_{jk})\in\text{\f C}'$ and
$$\multline \P(W^{0k}\in A_k,W^{0jk}\in A_{jk})=\P(W^{0k}_{F_y(\cdot)}\in F_y(A_k),W^{0jk}_{F_y(\cdot)}\in F_y(A_{jk}))\hfill\\ \hfill =\Pi_k\P(W^{0k}_{F_y(\cdot)}\in F_y(A_k))\times\Pi_{j<k}\P(W^{0jk}_{F_y(\cdot)}\in F_y(A_{jk}))\hfill\\ \hfill=\Pi_k\P(W^{0k}\in A_k)\times\Pi_{j<k}\P(W^{0jk}\in A_{jk}).\endmultline\tag3.20$$
Therefore the $W^{0k},W^{0jk}$ are independent and we have Theorem 1.3 in this case.

For $y>1$ we express the processes in the form of a matrix.  Let $W_n$ denote the $m\times m$ matrix with ${W_n}_{kk}=X_n^k$, and for $j<k$ ${W_n}_{jk}={W_n}_{kj}=Y_n^{jk}$.  Let $O_{n,1}$ and $O_1$ as in Theorem 3.1.  Let $\varphi_n(t)$  be as in Theorem 3.1 with $\varphi(t)=\min(t/(1-1/y)),1)$ as its a.s. limit.  Let $\psi_n(t)=\max(t,F_n(0))$ with $\psi(t)\equiv\max(t,1-1/y)$ as its a.s. limit. Let $B_m$ be the $m\times m$ matrix consisting of $1/\sqrt 2$'s on its diagonal and 1's on its off-diagonal elements.   Let $\underline X_m$ be the $m_1\times m$ matrix with $i$-th column $O_{n,1}^T\text{\bf x}_{n,i}$,  let $I_{m_1,s}$ be the $m_1\times m_1$ diagonal matrix consisting of 1's on its first $s$ diagonal entries, 0 on the remaining diagonal entries, and let $I_{m_1}$ be the $m_1\times m_1$ identity matrix.   Notice that $m_1=nF_n(0)$.  Denote ``$\circ$" as the Hadamard product.  Then we have
$$W_n(t)=\sqrt n B_m\circ\left(\underline X_m^TO_1I_{m_1,[m_1\varphi_n(t)]}O_1^T\underline X_m-\frac{[m_1\varphi_n(t)]}nI_{m_1}\right)-W_n(F_n(0))+
W_n(\psi_n(t)).$$
Let $W'$ be the weak limit of $W_n$ on a subsequence.  Then on this subsequence $W_n(\psi_n(\cdot))\darrow W'(\psi(\cdot))$, and $W_n(\psi_n(F_n(\cdot)))=W_n(F_n(\cdot))\darrow W'_{F_{y(\cdot)}}$, where the entries of $W'_{F_{y(\cdot)}}$ on and above the diagonal are independent copies of Brownian bridge, composed with $F_y$.  Confining to the interval $[1-1/y,1]$ these entries  will also be independent copies on $C[1-1/y,1]$.  If we 
define $F_y^{-1}$ just on $[1-1/y,1]$ we have for $X\in C'$ $X(F_y^{-1}(F_y))=X$.  Therefore from (3.19) we see that the entries on and above the diagonal of $W'_{F_y(\cdot)}$ are
independent.  For $X,Y\in C[1-1/y,1]$ $X\neq Y$ we have $X(F_y(\cdot))\neq Y(F_y(\cdot))$ so that the 1-1 condition of Theorem 8.3.7 of [6] is satisfied.   We also have $X(F_y(F^{-1}_y))=X$.  Therefore  we have from  (3.20) with the entries of $W'$ confined to $[1-1/y,1]$ and the sets Borel subsets of  $ C[1-1/y,1]$, the  entries of $W'$ on
$[1-1/y,1]$ on and above the diagonal are independent.  This uniquely determines the limiting
distribution, so we see that $W_n(\psi_n(\cdot))\darrow W^0(\psi(\cdot))$, where $W^0$ is Brownian bridge, with entries on and above the diagonal independent.

Let $\underline X_m=U_mR_m$ be the QR factorization of $\underline X_m$, where the columns of $U_m$ are orthonormal, and $R_m$ is $m\times m$ upper triangular, with nonnegative diagonal entries.  Extending (3.3) and (3.4) to all columns of $\underline X_m$ we have
$$R_m^TR_m=\underline X_m^T\underline X_m\iparrow(1-(1/y))I_m.$$
From this it is straightforward to prove
$$R_m\iparrow\sqrt{1-(1/y)}I_m.\tag3.21$$

Write
$$\sqrt n B_m\circ\left(\underline X_m^TO_1I_{m_1,[m_1\varphi_n(t)]}O_1^T\underline X_m-\frac{[m_1\varphi_n(t)]}nI_{m_1}\right)-W_n(F_n(0))\tag3.22$$
$$\multline =\frac1{\sqrt{F_n(0)}}B_m\circ R_m^T\sqrt{m_1}\left(U_m^TO_1I_{m_1,[m_1\varphi_n(t)]}O_1^TU_m-
\frac{[m_1\varphi_n(t)]}{m_1}I_{m_1}\right)R_m\hfill \\ \hfill+W_n(F_n(0))\left(\frac{[m_1\varphi_n(t)]}{m_1}-1\right).\endmultline$$
As in Theorem 2.1 of [16], we use Theorem 5.1 of [3] applied to 
$$\left(W_n,\sqrt{m_1}\left(U_m^TO_1I_{m_1,[m_1s]}O_1^TU_m-
\frac{[m_1s]}{m_1}I_{m_1}\right),R_m,F_n(0),\varphi_n,\psi_n\right).$$ 
We also apply Lemma 3.2 where $X'_n=(W_n,F_n(0))$, $\ell_n=m_1$, and $X''_{\ell_n}$ is the second component of the above six-tuple.  Therefore, from Theorem 1.2, (3.21), and (3.22) we
have 
$$W_n\darrow \sqrt{1-(1/y)}\hat W^0_{\varphi}+W^0_{1-(1/y)}(\varphi-1)+W^0_{\psi},$$
where $\hat W^0$ is an independent copy of $W^0$.   Since this limit is the same when $v_{11}$ is $N(0,1)$ we have this limit having independent elements on and above the diagonal.  This completes the proof of Theorem 1.3.

{\bf 4. Proof of Theorem 1.4} We first need the following:

Lemma 5 (Lemma 2.7 in Bai and Silverstein (1998).  For $X=(X_1,\ldots,X_n)^T$ i.i.d. standardized entries, and $C$, an $n\times n$
matrix, we have, for any $p\ge2$
$$\exp |X^*CX-\tr C|^p\leq K_p((\exp|X_1|^4\tr CC^*)^{p/2}+\exp|X_1|^{2p}\tr(CC^*)^{p/2}).$$

Suppose $C$, $n\times n$, is bounded in spectral norm and $X$ contains i.i.d. complex Gaussian entries.  Then for any $p\ge2$

$$\exp|X^*CX-\tr C|^p\leq K_p\|C\|^p((\exp|X_1|^4)^{p/2}n^{p/2}+\exp|X_1|^{2p}n)\leq K_pn^{p/2}.\tag4.1$$

Recalling $S_n=U_n\Lambda_nU_n^*$ in its spectral decomposition with eigenvalues arranged in nondecreasing order,  for any real $x$ let $\Lambda_n(x)$ denote the diagonal
matrix containing $nF_n(x)$ one's on the upper part of its diagonal.  Therefore $F_n(x)=(1/n)\tr\Lambda_n(x)$.  Notice that $G_n(x)=\bfv_n^*U_n\Lambda_n(x)U_n^*\bfv_n=\sum_{\lambda_k\leq x}|\bfu_k^*\bfv_n|^2$ where $U_n=(\bfu_1,\ldots,\bfu_n)$,
is the distribution function of a random variable which takes values $\lambda_1,\ldots,\lambda_n$ (eigenvalues of $S_n$) 
 with probabilities $|\bfu_1^*\bfv_n|^2,\ldots,|\bfu_n^*\bfv_n|^2$. Now, since $U_n^*\bfv_n$ is uniformly
distributed on the $n$ dimensional unit sphere in $\Bbb C^n$ it has the distribution of a normalized vector $\bfz_n$ of $n$ i.i.d. complex
Gaussian entries: $U^*\bfv_n\sim (1/\|\bfz_n\|)\bfz_n$. By (4.1) we have
$$\exp|(1/n)\bfz_n^*\Lambda_n(x)\bfz_n-F_n(x)|^4\leq Kn^{-2}.$$
Moreover
$$|G_n(x)-(1/n)\bfz_n^*\Lambda_n(x)\bfz_n|=(1/n)\bfz_n^*\Lambda_n(x)\bfz_n|n/\|\bfz_n\|^2-1|\asarrow0$$
by the strong law of large numbers.  Therefore we have with probability one, $G_n$ converges in
distribution to $F$, and the largest value in the support of $G_n$, namely the largest eigenvalue
of $S_n$, converges with probability one to $\lambda_{\max}$.  Therefore for any $\lambda>\lambda_{\max}$
with probability one, for all $n$ large $(\bfv_n^*(\lambda I-S_n)^{-1}\bfv_n,\bfv_n^*(\lambda I-S_n)^{-2}\bfv_n)$ exists and converges to 
$(\int(\lambda-x)^{-1}dF(x),(\int(\lambda-x)^{-2}dF(x))$.

Suppose that for all $\lambda>\lambda_{\max}$  $\int(\lambda-x)^{-1}dF(x)\leq1/\theta$.  Then necessarily
$\lim_{\lambda\to\lambda_{\max}^+}\int(\lambda-x)^{-1}dF(x)\leq1/\theta$, which means for all $\epsilon>0$
$\int(\lambda_{\max}+\epsilon-x)^{-1}dF(x)<1/\theta$  Since almost surely $\bfv_n^*((\lambda_{\max}+\epsilon)I-S_n)^{-1}\bfv_n\rightarrow\int(\lambda_{\max}+\epsilon-x)^{-1}dF(x)$,
we must have with probability one, for all $n$ large $\lambda_n^1<\lambda_{\max}+\epsilon$.  Since $\epsilon$ is arbitrary we must have
almost surely $\lambda_n^1\rightarrow\lambda_{\max}$.

Suppose now there exists $\lambda>\lambda_{\max}$ such that $\int(\lambda-x)^{-1}dF(x)>1/\theta$ Then let 
 $\lambda_1>\lambda_{\max}$ be the unique value such that $\int(\lambda_1-x)^{-1}dF(x)=1/\theta$.  For small $\epsilon>0$
$$\int(\lambda_1-\epsilon-x)^{-1}dF(x)>1/\theta\quad\text{and}\quad\int(\lambda_1+\epsilon-x)^{-1}dF(x)<1/\theta.$$
Since almost surely 
$$\multline \bfv_n^*((\lambda_1-\epsilon)I-S_n)^{-1}\bfv_n\rightarrow\int(\lambda_1-\epsilon-x)^{-1}dF(x)\\
\text{and}\quad
\bfv_n^*((\lambda_1+\epsilon)I-S_n)^{-1}\bfv_n\rightarrow\int(\lambda_1+\epsilon-x)^{-1}dF(x)\endmultline
$$ 
we have almost surely for all $n$ large $\lambda_1-\epsilon<\lambda_n^1<\lambda_1+\epsilon.$
Since $\epsilon$ is arbitrary we must have $\lambda_n^1\asarrow\lambda_1$.

For small $\epsilon>0$ we have with probability one, for all $n$ large
$$\bfv_n^*((\lambda_1+\epsilon)I-S_n)^{-2}\bfv_n\leq\bfv_n^*(\lambda_n^1I-S_n)^{-2}\bfv_n\leq\bfv_n^*((\lambda_1-\epsilon)I-S_n)^{-2}\bfv_n.$$
where the extremes approach almost surely $\int(\lambda_1+\epsilon-x)^{-2}dF(x)$, $\int(\lambda_1i\epsilon-x)^{-2}dF(x)$, respectively.  Since $\epsilon$ we have
$$\bfv_n^*(\lambda_n^1I-S_n)^{-2}\bfv_n\asarrow\int(\lambda_1-x)^{-2}dF(x),$$
which gives us (1.12).



Let $b\in(\lambda_{\max},\lambda_1)$ and $a=(\lambda_{\max}+b)/2$. 
Select $d>\lambda_1$.  Define for $t\in[b,d]$ $\Phi_n(t)\equiv b$, if  $\lambda_n^1\notin [b,d]$ and $\equiv\lambda_n^1$ if $\lambda_n^1\in [b,d]$.
Then $\Phi_n$ is a random element in $D_0[b,d]$, those elements of $D[b,d]$ whose range is also in $[b,d]$ and nondecreasing (pp. 144-145 of [3]).  Then with probability one, for all $n$ large, $\Phi_n\equiv\lambda_n^1$ and converges to 
$\lambda_1$.  

Identify $\bfv_n$ with $\bfx_{n,k'}$ in (1.4). Define $X_n^k(x)=X_n^{k,k'}(F_n(x))$.   We have $X_n^k$ a random element in 
$D[0,\infty)$, the set of all functions on $[0,\infty)$ having discontinuities of the first kind ([9]).   It is straightforward to 
extend the material in [3] pp. 144-145 and Theorem 4.4 to bounded nondecreasing functions in $D[0,\infty)$ to 
conclude that  $X_n^k(x)$ converges weakly to
$$W_{k,r}^0(F(x))+iW_{k,i}^0(F(x))\tag4.2$$
on $D_2[0,\infty)$ (two copies of $D[0,\infty)$) (Note: 
this is the only place where we need the limiting distribution function $F$ to be continuous).

 Let for $x\in [0,a]$
$$Y_n^k(x)=I_{\{\lambda_{\max}(S_n)\leq a\}}X_n^k(x),$$
where $I_A$ is the indicator function on the set $A$.  Then from Theorem 4.1 of [3] $Y_n^k$ converges weakly to (4.2) on $D_2[0,a]$ (two copies of $D[0,a]$).

Define the mapping $f$ from $D_2[0,a]$ to $C_2[b,d]$ (two copies of $C[b,d]$, the space of continuous functions on $[b,d]$) by 
$$f(X)=-\int_0^{a}(t-x)^{-2}X(x)dx\quad t\in[b,d].$$
Then
$$\multline f(Y_n^k)=-I_{\{\lambda_{\max}(S_n)\leq a\}}\int_0^{a}(t-x)^{-2}X_n^k(x)dx\hfill\\ \hfill I_{\{\lambda_{\max}(S_n)\leq a\}}\int_0^{a}(t-x)^{-1}dX_n^k(x)=I_{\{\lambda_{\max}(S_n)\leq a\}}\sqrt{2n}\bfx_{n,k}^*(tI-S_n)^{-1}\bfv_n.\endmultline$$
We claim that $f$ is a continuous mapping.  Suppose $X_n\to X$ in $D_2[0,a]$ in the Skorohod topology.  Then 
$X_n(s)\to X(s)$ for continuity points $s$ of $X$, and because $X$ lies in $D_2[0,a]$, this set is outside a set of Lebesgue measure 0.   Using the fact that convergence in the Skorohod topology renders the $X_n$ and $X$ uniformly bounded we have by the dominated convergence theorem
$$|f(X_n)-f(X)|\leq((b-\lambda_{\max})/2)^2\int_0^a|X_n(x)-X(x)|dx\to0,$$
uniformly for $t\in[b,d]$.  Therefore $f$ is continuous.

Therefore from Theorem 5.1 of [3] we have 
$$\multline I_{\{\lambda_{\max}(S_n)\leq a\}}\sqrt{2n}\bfx_{n,k}^*(tI-S_n)^{-1}\bfv_n\hfill\\ \hfill\darrow \int(t-x)^{-1}dW_{k,r}^0(F(x))+i \int(t-x)^{-1}dW_{k,i}^0(F(x))\endmultline$$
on $D_2[b,d]$.  
From the material on pp.144-145 of [3] we have
$$\multline I_{\{\lambda_{\max}(S_n)\leq a\}}\sqrt{2n}\bfx_{n,k}^*(\Phi_n I-S_n)^{-1}\bfv_n\hfill\\ \hfill\darrow \int(\lambda_1-x)^{-1}dW_{k,r}^0(F(x))+i \int(\lambda_1-x)^{-1}dW_{k,i}^0(F(x)).\endmultline$$
Using again Theorem 4.1 of [3] we get (1.10).
 
 We get the same result for $G_n$ in the real Gaussian case.  For the matrix $M_n$ it is proven in Section 3 that $G_n\darrow F_y$ i.p.  For the former the steps above follow identically, resulting in (1.13).   For the latter, since the finite result is distributional in nature we may as well assume $G_n\darrow F_y$ a.s. (since this is true on an appropriate subsequence  
 of an arbitrary subsequence of natural numbers). Thus we get (1.13) with $F=F_y$.
 \vskip.5in

\centerline{REFERENCES}
\medskip
\item{[1]} Baik, J., Ben Arous, G., and P\'ech\'e. (2005) Phase transition of the largest eigenvalue for non-null
complex sample covariance matrices. {\sl Ann. Probab. } {\bf33} 1643-1697.
\item{[2]} Baik, J, and Silverstein, J.W. (2006) Eigenvalues of large sample covariance matrices of spiked 
population models. {\sl J. Multivariate Anal. \bf97} 1382-1408
\item{[3]} Billingsley, P. (1968). {\sl Convergence of Probability
Measures}. Wiley, New York.
\item{[4]} Billingsley, P. (1995) {\sl Probability and Measure} Third Edition. Wiley,
New York.
\item{[5]}  Billingsley, P. (1999). {\sl Convergence of Probability
Measures} Second Edition. Wiley, New York.
\item{[6]} Cohn, D.L. (1980) {\sl Measure Theory} Birkhauser Boston.
\item{[7]} Grenander, U. and Silverstein, J.W. (1977). Spectral
analysis of networks with random topologies. {\sl SIAM J. Appl. Math.
\bf37} 499-519.
\item{[8]} Jonsson, D. (1982). Some limit theorems for the
eigenvalues of a sample covariance matrix. {\sl J. Multivariate Anal.
\bf12} 1-38.
\item{[9]} Lindvall, T. (1973). Weak convergence of probability
measures and random functions in the function space $D[0,\infty)$).
{\sl J. Appl. Probab. \bf10} 109-121.
\item{[10]} Mar\v cenko, V.A. and Pastur, L.A. (1967). Distribution of eigenvalues for
some sets of random matrices. {\sl Math. USSR-Sb.} {\bf1}, 457-483.
\item{[11]} Rao, R. and Silverstein, J.W. (2010) Fundamental limit of sample generalized eigenvalue based
detection of signals in noise using relatively few signal-bearing and noise-only samples 
 {\sl IEEE Journal of Selected Topics in Signal Processing} {\bf3} 468-480.
\item{[12]} Silverstein, J.W. (1979). On the randomnes of eigenvectors
generated from networks with random topologies. {\sl SIAM J. Appl.
Math. \bf37} 235-245.
\item{[13]} Silverstein, J.W. (1981). Describing the behavior of
random matrices using sequences of measures on orthogonal groups.
 {\sl SIAM J. Math. Anal. \bf12} 274-281.
\item{[14]} Silverstein, J.W. (1984). Some limit theorems on the
eigenvectors of large dimensional sample covariance matrices. 
 {\sl J. Multivariate Anal.} {\bf15} 295-324.
\item{[15]} Silverstein, J.W. (1989). On the eigenvectors of large
dimensional sample covariance matrices. {\sl J. Multivariate Anal.}
{\bf 30} 1-16.
\item{[16]} Silverstein, J.W. (1990). Weak convergence of random functions defined by the
eigenvectors of sample covariance matrices.  {\sl Ann. Probab. \bf18} 1174-1194.
\item{[17]} Silverstein, J.W. (1995). Strong convergence of the empirical distribution of eigenvalues of large
dimensional random matrices  {\sl J. Multivariate Anal. \bf55} 331-339.
\item{[18]} Wachter, K.W. (1978). The strong limits of random matrix
spectra for sample matrices of independent elements. {\sl Ann.
Probab. \bf6} 1-18.
\item{[19]} Yin, Y.Q. (1986). Limiting spectral distribution for a
class of random matrices. {\sl J. Multivariate Anal. \bf20} 50-68.
\item{[20]} Yin, Y.Q., Bai, Z.D., and Krishnaiah, P.R. (1988). On limit of
the largest eigenvalue of the large dimensional sample covariance
matrix. {\sl Probab. Th. Rel. Fields \bf78} 509-521.

\bye

Define for $j\leq m$ and $t\in[b,d]$
$$X_n^{j}(t)=I_{\{\lambda_{\max}(S_n)\leq(\lambda_{\max}+b)/2\}} \sqrt{2n}x_j^*(tI-S_n)^{-1}u_{n}.$$
$X_n^j$ can be viewed as a two-dimensional continuous process on $[b,d]$.    Since with probability one $\sqrt{2n}x_j^*(\lambda_n^1I-S_n)^{-1}u_n=X_n^j(\Phi(t))$ for all $n$ large, it follows from the
material on pp. 144-145 of Billingsley (1968) that (14) is proven once it is shown that the real and imaginary parts of $X_n^j$ converge weakly in $C[b,d]$ to
$$\int_{\lambda_{\min}}^{\lambda_{\max}}(t-x)^{-1}dW_{j,r}^0(F(x))\quad\text{and}\quad \int_{\lambda_{\min}}^{\lambda_{\max}}(t-x)^{-1}dW_{j,i}^0(F(x)),$$
Using Theorem 5.5 of Billingsley (1968) on $f_n^t(x)=I_{\{\lambda_{\max}(S_n)\leq(\lambda_{\max}+b)/2\}}(t-x)^{-1}$, $t\in[b,d]$, and the results of the previous processes to Brownian bridge, we have the convergence
of the finite dimensional distributions to the proper limits. To show tightness we verify (12.51) of Billingsley (1968).  Since
$$\max((\Re X_n^j(t_2)-\Re X_n^j(t_1))^2,(\Im X_n^j(t_2)-\Im X_n^j(t_1))^2)\leq |X_n^j(t_2)-X_n^j(t_1)|^2,$$
we will work with the latter.

Let $\lambda_1\leq\dots\leq\lambda_n$ be the eigenvalues of $S_n$ and $\Lambda_n$ the diagonal matrix containing the eigenvalues.  Write 
$$\multline X_n^j(t)=I_{\{\lambda_{\max}(S_n)\leq(\lambda_{\max}+b)/2\}} \sqrt{2n}u^*(tI-\Lambda_n)^{-1}v\hfill\\ \hfill=I_{\{\lambda_{\max}(S_n)\leq(\lambda_{\max}+b)/2\}} \sqrt{2n}\sum_i\bar u_iv_i(t-\lambda_i)^{-1},
\endmultline$$
where $u$ and $v$ can be viewed as the first two columns of a Haar distributed matrix.  From symmetry we have
$$0=\exp \bar u_iv_i(\sum_ku_k\bar v_k)=\exp(|u_i|^2|v_i|^2)+\sum_{k\neq i}\exp(\bar u_iv_iu_k\bar v_k)=\exp(|u_i|^2|v_i|^2)+(n-1)\exp(\bar u_iv_iu_k\bar v_k)$$
for any $i\neq k$.  Also from symmetry $\exp(|u_i|^2|v_i|^2)=\exp(|u_1|^2|v_1|^2)$ and $\exp(\bar u_iv_iu_k\bar v_k)=\exp(\bar u_1v_1u_2\bar v_2)$.
Therefore we have for $t_1,t_2\in [b,d]$
$$\multline
\exp|X_n^j(t_2)-X_n^j(t_1)|^2=(t_2-t_1)^22n\exp(|u_1|^2|v_1|^2)\exp I_{\{\lambda_{\max}(S_n)\leq(\lambda_{\max}+b)/2\}}\left(\sum_i(t_2-\lambda_i)^{-2}(t_1-\lambda_i)^{-2}\right.\hfill\\ \hfill\left.-(n-1)^{-1}\sum_{i\neq k}(t_2-\lambda_i)^{-1}(t_1-\lambda_i)^{-1}(t_2-\lambda_k)^{-1}(t_1-\lambda_k)^{-1})\right)\endmultline$$
$$\leq Kn^2(t_2-t_1)^2\exp(|u_1|^2|v_1|^2).$$
Now $u_1$ and $v_1$ can also be considered as the first two entries in a row of a Haar distributed matrix.  Therefore, again from symmetry
$$n^{-1}=\exp|u_1|^2=\exp|u_1|^4+(n-1)\exp(|u_1|^2|v_1|^2).$$
Therefore we have
$$\exp|X_n^j(t_2)-X_n^j(t_1)|^2\leq K(t_2-t_1)^2,$$
so that (12.51) holds for the real and imaginary parts of $X_n^j$ and we are done.

In the real case we will not have $\sqrt 2$ in (14).  The results are identical.
\bye

 \bye
 
The main result in this paper is the following.  For each $n$, let $\bfx_{n,1},\ldots,\bfx_{n,m} $ denote orthogonal nonrandom unit vectors in ${\Bbb R}^n$ and let 
$\text{\bf y}_k=(y_{k,1},\ldots,y_{k,n})^T=O^T\text{\bf x}_{n,k}$, $k=1,\ldots,m$.  For each $k$ define $X_n^k$, a random element in $D[0,1]$ to be (1.1) with $y_i$ replaced 
with $y_{k,i}$  For $1\leq j<k\leq m$ define $Y_n^{jk}$, a random element of $D[0,1]$, to be
$$Y_n^{jk}(t)=\sqrt{n}\sum_{i=1}^{[nt]}y_{j,i}y_{k,i},$$
Then the random functions $X_n^k,Y_n^{jk}$, $1\leq j<k\leq m$ converge weakly in $D_d^1$, $d=m(m+1)/2$, to independent Brownian Bridges. 

Let for $d\ge2$ an integer, and $b\ge1$ $D_d^b=\Pi_{i=1}^dD[0,b]$, and $\text{\f T}_d^b$ denote the smallest $\sigma$-field on $D_d^b$ in which convergence of elements in $D_d^b$ is equivalent to component-wise convergence.

We will prove

{\smc Theorem 1.2} Assume $v_{11}$ is symmetrically distributed about 0, $\exp v_{11}^4<\infty$, and each $\bfx_{n,k}=(\pm1/\sqrt n,\ldots,\pm1/\sqrt n)^T$.   Then the random functions $X_n^k,Y_n^{jk}$, $1\leq j<k\leq m$ converge weakly in $D_d^1$, $d=m(m+1)/2$, to independent Brownian Bridges. 

 {\bf 1. Introduction.}	 Let $\{v_{ij}\}$, $i,j=1,2,\ldots$,
be i.i.d. random variables with \text{\pe E}$(v_{11})=0$, 
and for each $n$ let
${M_n=\frac1sV_nV_n^T}$, where $V_n=(v_{ij})$, $i=1,2,\ldots,n$,
$j=1,2,\ldots,s=s(n)$, and $n/s\rightarrow y>0$ as
$n\rightarrow\infty$. The symmetric, nonnegative definite matrix $M_n$ can 
be viewed as the sample	covariance matrix of $s$ samples of an $n$
dimensional random vector having i.i.d. components distributed the
same as $v_{11}$ (assuming knowledge of the common mean being zero). 
The spectral behavior of $M_n$ for $n$ large is important to areas of
multivariate analysis (including principal component analysis,
regression, and signal processing) where $n$ and $s$ are the same
order of magnitude, so that standard asymptotic analysis cannot be
applied.  Although eigenvalue results will be discussed,  this paper is 
chiefly concerned with the behavior of the eigenvectors of $M_n$.  
It continues the analysis begun in [8].

Throughout the following, $O_n\Lambda_n O_n^T$ will denote the spectral 
decomposition of $M_n$, where the eigenvalues of $M_n$ are arranged
along the diagonal of $\Lambda_n$ in nondecreasing order.  The orthogonal 
matrix $O_n$, the columns being eigenvectors of $M_n$, is not uniquely
determined, owing to the multiplicities of the eigenvalues and the
direction any eigenvector can take on. 	This problem will be
addressed later on.  For now it is sufficient to mention that, by
appropriately enlarging the sample space where the $v_{ij}$'s are
defined, it is possible to construct $O_n$ measurable in a natural
manner from the eigenspaces associated with $M_n$. 

Results previously obtained suggest similarity of behavior of the 
eigenvectors of $M_n$ for large	$n$ to the eigenvectors of matrices of 
Wishart type, that is, when $v_{11}$ is normally distributed 
([8],[9],[10],[11]). In this case it is well-known that $O_n$
induces the Haar (uniform) measure on ${\Cal O}_n$, the $n\times n$
orthogonal group.  In [9] it is conjectured that, for general
$v_{11}$ and for $n$ large, $O_n$ is somehow close to being Haar
distributed.  The attempt to make the notion of
closeness more precise has led in [9] to an investigation into the behavior	
of random elements, $X_n$, of $D[0,1]$ (the space of $r.c.l.l.$ functions on
[0,1]) defined by the eigenvectors.  They are constructed as follows:

For each $n$ let $\vec x_n\in R^n$, $\|x_n\|=1$, be nonrandom, and let 
$\vec y_n=(y_1,y_2,\ldots,y_n)^T=O_n^T\vec x_n$.  Then, for $t\in [0,1]$
$$X_n(t)\equiv\sqrt{\tfrac n2}
\sum_{i=1}^{[nt]}(y_i^2-\tfrac1n)\quad\quad([a]\equiv\text{greatest
integer}\leq a).$$

The importance of $X_n$ to understanding the behavior of the
eigenvectors of $M_n$ stems from three facts.  First, the behavior of
$X_n$ for all $\vec x_n$ reflects to some degree the uniformity or
nonuniformity of $O_n$.  If $O_n$ were Haar distributed, then $\vec y_n$ 
would be uniformly distributed over the unit sphere in $R^n$, rendering an 
identifiable distribution for $X_n$, invariant across unit $\vec x_n\in R^n$.
On the other hand, significant departure from Haar measure would be
suspected if the distribution of $X_n$ depended strongly on $\vec
x_n$.   Second, we are able to compare some aspects
of the distribution of $O_n$ for all $n$ on a common space, namely
$D[0,1]$.  Third, the limiting behavior of $X_n$ is known when $O_n$
is Haar distributed.  Indeed, in this case the distribution  of $\vec
y_n$, being uniform on the unit sphere in $R^n$, is  
the same as that of a normalized vector of i.i.d. mean-zero
Gaussian components.  Upon applying standard results on weak
convergence of measures, it is straightforward to show
$$X_n\darrow W^\circ\quad\text{ as }
n\rightarrow\infty\leqno(1.1)$$
(\cal D denoting weak convergence in $D[0,1]$) where $W^\circ$ is
Brownian bridge ([9]).  Thus, for arbitrary $v_{11}$, (1.1) holding for 
all $\{\vec x_n\}$, $\|x_n\|=1$, can be viewed as evidence supporting
the conjecture that $O_n$ is close to being uniformly distributed in 
${\Cal O}_n$ for $n$ large.

 Verifying (1.1) more generally is difficult since there is very little 
useful direct information available on the variables $y_1^2,\ldots,y_n^2$.  
However, under the assumption \text{\pe E}$(v_{11}^4)<\infty$, the results in 
[9],[10],[11] reduce the problem to verifying \underbar{tightness} of 
$\{X_n\}$.  For the following we may, without loss of generality, 
assume \text{\pe E}$(v_{11}^2)=1$.  The results
are limit theorems on random variables defined by $\{M_n\}$.  From
one of the theorems ([11], to be given below) it follows that for 
\text{\pe E}$(v_{11}^4)=3$, any weakly convergent subsequence of $\{X_n\}$
converges to $W^\circ$ for any sequence $\{\vec x_n\}$ of unit vectors, 
while if \text{\pe E}$(v_{11}^4)\neq 3$, there exists sequences 
$\{\vec x_n\}$ for which $\{X_n\}$ fails to converge weakly.  This
suggests some further similarity of $v_{11}$ to $N(0,1)$ may be
necessary. It is remarked here that the theorem also implies, for
finite \text{\pe E}$(v_{11}^4)$, the necessity of \text{\pe E}$(v_{11})=0$
in order for (1.1) to hold for $\vec x_n=(1,0,\ldots,0)^T$.

The main purpose of this paper is to establish the following partial solution
to the problem:
\medskip

{\smc Theorem 1.1}.  {\sl Assume $v_{11}$ is symmetric (that is, symmetrically
distributed about 0), and \text{\pe E}$(v_{11}^4)<\infty$.  Then
(1.1) holds for $\vec x_n=(\pm\frac1{\sqrt n},\pm\frac1{\sqrt n},\ldots,
\pm\frac1{\sqrt n})^T$.}
\medskip

 From the theorem one can easily argue other choices of $\vec x_n$ for
which (1.1) holds, namely vectors close enough to those in the
theorem so that the resulting $X_n$ approaches in the
Skorohod metric random functions satisfying (1.1).
It will become apparent that the techniques used in the proof of
Theorem 1.1 cannot easily be extended to $\vec x_n$ having more
variability in the magnitude of	its components, while the symmetry 
requirement may be weakened with a deeper analysis.  At present the 
possibility exists that only for $v_{11}$ mean-zero Gaussian will (1.1) 
be satisfied for all $\{\vec x_n\}$.
 
However, from Theorem 1.1 and the previously mentioned results emerges the 
possibility of classifying the distribution of $O_n$ into varying degrees 
of closeness to Haar measure.  The 
eigenvectors of $M_n$ with $v_{11}$ symmetric and fourth moment	finite
display a certain amount of uniform behavior, and $O_n$
can possibly be even more closely related to Haar measure if 
{\pe E}$(v_{11}^4)/[\text{\pe E}(v_{11}^2)]^2=3$. 
As will be seen below when it is
formally stated, the limit theorem in [11] itself can be viewed as 
demonstrating varying degrees of similarity to Haar measure. 

The proof of Theorem 1.1 relies on two results on the eigenvalues of
$M_n$ and on a modification of the limit theorem in [11].  Let $F_n$ denote the
empirical distribution function of the eigenvalues of $M_n$ (that is,
$F_n(x)=(1/n)\times\text{(number of eigenvalues of }M_n\leq x)$,
where we may as well assume $x\geq0$).  If {\pe Var}$(v_{11})=1$ (no
other assumption on the moments), then it is known ([4],[5],[12],[13])
that, for every $x\geq0$, 
$$F_n(x)\asarrow F_y(x)\quad\text{ as }n\rightarrow\infty,\leqno(1.2)$$ 
where $F_y$ is a continuous,
nonrandom probability distribution function depending only on $y$,
having a density with support on $[(1-\sqrt y)^2,(1+\sqrt y)^2]$, and
for $y>1$, $F_y$ places mass $1-1/y$ at 0.  Moreover, if {\pe E}$(v_{11})=0$, 
{\pe E}$(v_{11}^2)=1$, then $\lambda_{\text{max}}(M_n)$, the largest
eigenvalue of $M_n$, satisfies
$$\lambda_{\text{max}}(M_n)\asarrow(1+\sqrt y)^2\quad\text{ as }
n\rightarrow\infty\leqno(1.3)$$
if and only if {\pe E}$(v_{11}^4)<\infty$ ([1],[3],[14]).

For the theorem in [11], we first make the following
observation.  It is straightforward to show that (1.1),(1.2),(1.3)
imply $$X_n(F_n(x))\darrow W^y_x\equiv W^\circ(F_y(x))\leqno(1.4)$$
on $D[0,\infty)$.  The proof of Theorem 1.1 essentially verifies the
truth of the implication in the other direction and then the truth of
(1.4). An extension of the theorem in [11] is needed for the latter to 
establish the uniqueness of any weakly convergent subsequence.   
The theorem states:
$$\left\{\,\sqrt{\tfrac{n}{2}}(\vec x_n^TM_n^r\vec x_n - 
\tfrac{1}{n}\text{tr}(M_n^r))\,\right\}_{r=1}^\infty\,\,=\,\,\left\{\,
\int_0^\infty x^r\,dX_n(F_n(x))\right\}_{r=1}^\infty\leqno(1.5)$$ 
$$\left\{\,-\int_0^\infty rx^{r-1}\,X_n(F_n(x))dx\right\}_{r=1}^\infty\darrow
\left\{\int_{(1-\sqrt{y})^2}^{(1+\sqrt{y})^2}x^r\,dW_x^y\right\}_{r=1}^\infty 
\qquad \text{as $n\rightarrow\infty$}$$
($\Cal D$ denoting weak convergence on $R^\infty$) for every sequence
$\{\vec x_n\}$ of unit vectors if and only if {\pe E}$(v_{11})=0$, 
{\pe E}$(v_{11}^2)=1$, and {\pe E}$(v_{11}^4)=3$ (we remark here that
the limiting random variables in (1.5) are well defined stochastic
integrals, being jointly normal each with mean 0).  The proof of this
theorem will be modified to show (1.5) still holds under the assumptions 
of Theorem 1.1 (without a condition on the fourth moment of
$v_{11}$ other than it being finite).

The proof will be carried out in the next three
sections.  Section 2 presents a formal description of $O_n$
to account for the ambiguities mentioned earlier, followed by a
result which converts the problem to one of showing weak
convergence of $X_n(F_n(\cdot))$ on $D[0,\infty)$.  Section 3 contains results 
on random elements in $D[0,b]$ for any $b>0$, which are extensions of certain 
criteria for weak convergence given in [2].  In section 4 the proof is 
completed by showing the conditions in section 3 are met.  Some of the
results will be stated more generally than presently needed to render them 
applicable for future use.
\vskip.1in
 \def\sqr#1#2{{\vcenter{\vbox{\hrule height.#2pt
\hbox{\vrule width.#2pt height#1pt \kern#1pt
\vrule width.#2pt}\hrule height.#2pt}}}}
\def\square{\mathchoice\sqr54\sqr54\sqr{2.1}3\sqr{1.5}3}

{\bf 2. Converting to $D[0,\infty)$. } Let us first give a more 
detailed description of the distribution of
$O_n$ which will lead us to a concrete construction of $\vec
y_n$.  For an eigenvalue $\lambda$ with multiplicity $r$ we assume the
corresponding $r$ columns of $O_n$ to be generated uniformly, that is, 
its distribution is the same as $O_{n,r}O_r$ where $O_{n,r}$ is
$n\times r$ containing $r$ orthonormal columns from the eigenspace of
$\lambda$, and $O_r\in{\cal O}_r$ is Haar distributed, independent of
$M_n$.  The $O_r$'s corresponding to distinct eigenvalues are also
assumed to be independent.  The coordinates of $\vec y_n$ corresponding
to $\lambda$ are then of the form $$(O_{n,r}O_r)^T\vec
x_n=O_r^TO_{n,r}^T\vec x_n=\|O_{n,r}^T\vec x_n\|\vec w_r$$ where $\vec
w_r$ is uniformly distributed on the unit sphere in $R^r$.  We
will use the fact that the distribution of $\vec w_r$ is the same as
that of a normalized vector of i.i.d. mean zero Gaussian components. 
Notice that $\|O_{n,r}^T\vec x_n\|$ is the length of the projection
of $x_n$ on the eigenspace of $\lambda$.

Thus, $\vec y_n$ can be represented as follows:

Enlarge the sample space defining $M_n$ to allow the construction of
$z_1,z_2,\ldots,z_n$, i.i.d. N(0,1) random variables independent of
$M_n$.  For a given $M_n$ let
$\lambda_{(1)}<\lambda_{(2)}<\cdots<\lambda_{(t)}$ be the $t$
distinct eigenvalues with multiplicities $m_1,m_2,\ldots,m_t$.  For
$i=1,2,\ldots,t$ let $a_i$ be the length of the projection of $\vec
x_n$ on the eigenspace of $\lambda_{(i)}$.  Define $m_0=0$.  Then, for each
$i$ we define the coordinates
$$(y_{m_1+\cdots+m_{i-1}+1},y_{m_1+\cdots+m_{i-1}+2},\ldots,
y_{m_1+\cdots+m_i})$$ of $\vec y_n$ to be the respective coordinates
of
$${a_i(z_{m_1+\cdots+m_{i-1}+1},z_{m_1+\cdots+m_{i-1}+2},\ldots,
z_{m_1+\cdots+m_i})}\left/{\sqrt{\dsize
\sum_{k=1}^{m_i}
z_{m_1+\cdots+m_{i-1}+k}^2}}\right..\leqno(2.1)$$

We are now in a position to prove 

{\smc Theorem 2.1}.  {\sl $X_n(F_n(\cdot))\darrow W_{F_y(\cdot)}^{\circ}$ in
$D[0,\infty)$, $F_n(x)\parrow F_y(x)$, and $\lambda_{\max}\parrow
(1+\sqrt y)^2$ $\Longrightarrow$ $X_n\darrow W^{\circ}$.}

{\smc Proof}.  We assume the reader is familiar with the basic results 
in [2] of showing weak convergence of random elements of a metric space
(most notably Theorems 4.1, 4.4, and Corollary 1 to Theorem 5.1), in
particular, the results on the function spaces $D[0,1]$ and $C[0,1]$.  For the
topology and conditions of weak convergence in $D[0,\infty)$ see [6]. 
For our purposes, the only information needed regarding $D[0,\infty)$
beyond that of [2] is the fact that weak convergence of a sequence of
random functions on $D[0,\infty)$ is equivalent to the following: for
every $B>0$ there exists a $b>B$ such that the sequence on $D[0,b]$
(under the natural projection) converges weakly.  Let $\rho$
denote the sup metric used on $C[0,1]$ and $D[0,1]$ (used only in the
latter when limiting distributions lie in $C[0,1]$ with probability
1), that is, for $x,y\in D[0,1]$ $$\rho(x,y)=\sup_{t\in
[0,1]}|x(t)-y(t)|.$$

We need one further general result on weak convergence, which is an
extension of the material on pp. 144-145 in [2] concerning random
changes of time.  Let $$\underline D[0,1]=\{x\in D[0,1]:\,x \text{ is
nonnegative and nondecreasing}\}$$ Since it is a
closed subset of $D[0,1]$ we take the topology of $\underline D[0,1]$
to be the Skorohod topology of $D[0,1]$ relativized to it.  The
mapping $$h:D[0,\infty)\times\underline D[0,1]\longrightarrow
D[0,1]$$ defined by $h(x,\varphi)=x\circ\varphi$ is measurable (same
argument as in [2], p. 232, except $i$ in (39) now ranges on all
natural numbers).  It is a simple matter to show that $h$ is
continuous for each $$(x,\varphi)\in C[0,\infty)\times C[0,1]\cap
\underline D[0,1].$$  Therefore, we have (by Corollary 1 to Theorem
5.1 of [2]) 
$$(Y_n,\Phi_n)\darrow (Y,\Phi)\text{ in }
D[0,\infty)\times\underline D[0,1],\,\text{\pe P}
\bigl(Y \in C[0,\infty)\bigr)
=\text{\pe P}\bigl(\Phi\in C[0,1]\bigr)=1\leqno(2.2)$$
$$\Longrightarrow Y_n\circ\Phi_n\darrow Y\circ\Phi\quad\text{ in }D[0,1].$$  

We can now proceed with the proof of the theorem.  Since we are
only concerned with distributional results we may as well assume that 
for all $x\geq 0$, $F_n(x)\asarrow F_y(x)$ and $\lambda_{\max}\asarrow
(1+\sqrt y)^2$ (for this is true on an appropriate subsequence of an
arbitrary subsequence of the natural numbers).  For $t\in [0,1]$
let $F_n^{-1}(t)$ = largest $\lambda_j$ such
that $F_n(\lambda_j)\leq t$ (0 for $t < F_n(0)$).  We have
$X_n\bigl(F_n(F_n^{-1}(t))\bigr)=X_n(t)$ except on intervals 
$[\frac mn,\frac{m+1}n)$ where $\lambda_m=\lambda_{m+1}$.  Let
$F_y^{-1}(t)$ be the inverse of $F_y(x)$ for ${x\in\bigl((1-\sqrt
y)^2,(1+\sqrt y)^2\bigr]}$

We consider first the case $y\leq 1$.  Let $F_y^{-1}(0)=(1-\sqrt
y)^2$.  It is straightforward to show for all $t\in (0,1]$,
 $F_n^{-1}(t)\asarrow F_y^{-1}(t)$.  Let
$\tilde{F}_n^{-1}(t) = \max\bigl( (1-\sqrt y)^2,F_n^{-1}(t)\bigr)$. 
Then, for all $t\in[0,1]$, $\tilde F_n^{-1}(t)\asarrow F_y^{-1}(t)$, and
since $\lambda_n\asarrow(1+\sqrt y)^2$ we have 
$\rho\bigl(\tilde{F}_n^{-1},F_y^{-1}\bigr)\asarrow 0$.
Therefore, from (2.2) (and Theorem 4.4 of [2]) we have 
$$X_n\bigl(F_n(\tilde{F}_n^{-1}(\cdot))\bigr)\darrow
W_{F_y(F_y^{-1}(\cdot))}^{\circ}=W^{\circ}\quad\text{in $D[0,1]$}.$$

Since $F_y(x)=0$ for $x\in[0,(1-\sqrt y)^2]$ we have 
$X_n(F_n(\cdot))\darrow 0$ in $D[0,(1-\sqrt y)^2]$, which implies
$X_n(F_n(\cdot))\parrow 0$ in $D[0,(1-\sqrt y)^2]$, and since the zero
function lies in $C[0,(1-\sqrt y)^2]$ we conclude that
$$\sup_{t\in[0,(1-\sqrt y)^2]}|X_n(F_n(x))|\parrow 0.$$
We have then $$\rho\bigl(X_n\bigl(F_n(F_n^{-1}(\cdot))\bigr),
X_n\bigl(F_n(\tilde{F}_n^{-1}(\cdot))\bigr)\bigr)\leq 
2\times\sup_{t\in[0,(1-\sqrt y)^2]}|X_n(F_n(x))|\parrow 0.$$
Therefore, we have (by Theorem 4.1 of [2])
$$X_n\bigl(F_n(F_n^{-1}(\cdot))\bigr)\darrow W^{\circ}\quad\text{in
$D[0,1]$}.$$

Notice if $v_{11}$ has a density then we would be done with this
case of the proof since for $n\leq s$ the eigenvalues would be
distinct with probability 1, 
so that $X_n\bigl(F_n(F_n^{-1}(\cdot))\bigr)=X_n(\cdot)$ almost
surely. However, for more general $v_{11}$, the multiplicities
of the eigenvalues need to be accounted for.

For each $M_n$ let
$\lambda_{(1)}<\lambda_{(2)}<\cdots<\lambda_{(t)}$,
$(m_1,m_2,\ldots,m_t)$, and $(a_1,a_2,\ldots,a_t)$ be defined as
above.  We have from (2.1)
$$\rho\bigl(X_n(\cdot),X_n(F_n(F_n^{-1}(\cdot)))\bigr)=
\max\Sb 1\leq i\leq t\\1\leq j\leq m_i\endSb\sqrt{\frac n2}\left|
a_i^2\frac{\dsize\sum_{\ell=1}^jz_{m_1+\cdots+m_{i-1}+\ell}^2}
{\dsize\sum_{k=1}^{m_i}z_{m_1+\cdots+m_{i-1}+k}^2}-\frac jn
\right|.\leqno(2.3)$$

The measurable function $h$ on $D[0,1]$ defined by
$$h(x)=\rho(x(t),x(t-0))$$ is continuous on $C[0,1]$ (note that $h(x)
={\dsize\lim_{\delta\downarrow 0}}w(x,\delta)$ where $w(x,\delta)$ is
the modulus of continuity of $x$) and is identically zero on $C[0,1]$.
Therefore (using Corollary 1 to Theorem 5.1 of [2])
$h\bigl(X_n(F_n(F_n^{-1}(\cdot)))\bigr)\darrow 0$, which is
equivalent to 
$$\max_{1\leq i\leq t}\sqrt{\frac n2}\left|a_i^2-\frac{m_i}n\right|\parrow
0.\leqno(2.4)$$

For each $i\leq t$ and $j\leq m_i$ we have 
$$\sqrt{\frac n2}\left(
a_i^2\frac{\dsize\sum_{\ell=1}^jz_{m_1+\cdots+m_{i-1}+\ell}^2}
{\dsize\sum_{k=1}^{m_i}z_{m_1+\cdots+m_{i-1}+k}^2}-\frac jn
\right)\quad=$$
$$\sqrt{\frac n2}(a_i^2-\frac {m_i}n)
\frac{\dsize\sum_{\ell=1}^jz_{m_1+\cdots+m_{i-1}+\ell}^2}
{\dsize\sum_{k=1}^{m_i}z_{m_1+\cdots+m_{i-1}+k}^2}\quad\quad+\leqno(a)$$
$$\sqrt{\frac n2}\frac {m_i}n\left(
a_i^2\frac{\dsize\sum_{\ell=1}^jz_{m_1+\cdots+m_{i-1}+\ell}^2}
{\dsize\sum_{k=1}^{m_i}z_{m_1+\cdots+m_{i-1}+k}^2}-\frac j{m_i}
\right).\leqno(b)$$

From (2.4) we have the maximum of the absolute value of (a) over
$1\leq i\leq t$ converges in probability to zero.  For the maximum of
(b) we see that the ratio of chi-square random variables is beta
distributed with parameters $p=j/2$, $q=(m_i-j)/2$.  Such a
random variable with $p=r/2$, $q=(m-r)/2$ has mean 
$r/m$ and fourth central moment bounded by $Cr^2/m^4$
where $C$ does not depend on $r$ and $m$.  Let $b_{m_i,j}$ represent
the expression in parentheses in (b).  Let $\epsilon>0$ be arbitrary.
We use Theorem 12.2 of [2] after making the following associations:
$S_j=\sqrt{m_i}b_{m_i,j}$, $m=m_i$, $u_{\ell}=\sqrt C/m_i$,
 $\gamma=4$, $\alpha=2$, and $\lambda=\epsilon\sqrt{2n/m_i}$.
We then have
$$\text{\pe P}\left(\max_{1\leq j\leq m_i}\left|\sqrt{\frac n2}\frac{m_i}nb_{m_i,j}
\right|>\epsilon\,\,\,\biggl|M_n\right)\leq\frac{C'm_i^2}
{4n^2\epsilon^4}.$$
By Boole's inequality we have
$$\text{\pe P}\left(\max\Sb 1\leq i\leq t\\1\leq j\leq m_i\endSb
\left|\sqrt{\frac n2}\frac{m_i}nb_{m_i,j}
\right|>\epsilon\,\,\,\biggl|M_n\right)\leq\frac{C'}
{4\epsilon^4}\max_{1\leq i\leq t}\frac{m_i}n.$$
Therefore
$$\text{\pe P}\left(\max\Sb 1\leq i\leq t\\1\leq j\leq m_i\endSb
{\left|\sqrt{\frac n2}\frac{m_i}nb_{m_i,j}
\right|}>\epsilon\right)\leq\frac{C'}
{4\epsilon^4}\text{\pe E}\left(\max_{1\leq i\leq t}\frac{m_i}n\right).\leqno(2.5)$$

We have $F_n(x)\asarrow F_y(x)\,\,\Longrightarrow
{\dsize\sup_{x\in[0,\infty)}}|F_n(x)-F_y(x)|\asarrow 0
\,\,\Longrightarrow$ (since $F_y$ is continuous on $(-\infty,\infty)$)
 ${\dsize\sup_{x\in[0,\infty)}}|F_n(x)-F_n(x-0)|\asarrow
0$, which is equivalent to ${\dsize\max_{1\leq i\leq
t}}{m_i}/n\asarrow 0$.  Therefore, by the dominated convergence
theorem, we have the left hand side of $(2.5)\longrightarrow 0$.  
We therefore have $(2.3)\parrow 0$ and we conclude (again from Theorem 4.1 
of [2]) that $X_n\darrow W^{\circ}$ in $D[0,1]$.

For $y > 1$ we assume $n$ is sufficiently large so that $n/s>1$.  Then
$F_n(0)=m_1/n\ge1-(s/n)>0$.    For
${t\in [0,1-(1/y)]}$ define $F_y^{-1}(t) = (1-\sqrt y )^2$.  For  $t\in
(1-(1/y),1]$ we have $F_n^{-1}(t)\asarrow F_y^{-1}(t)$.  Define as
before
$\tilde{F}_n^{-1}(t) = \max\bigl( (1-\sqrt y)^2,F_n^{-1}(t)\bigr)$. 
Again, $\rho(\tilde{F}_n^{-1},F_y^{-1})\asarrow 0$, and from (2.2)
(and Theorem 4.4 of [2]) we have
$$X_n\bigl(F_n(\tilde{F}_n^{-1}(t))\bigr)\darrow
W_{F_y(F_y^{-1}(t))}^\circ\,=\cases W_{1-(1/y)}^{\circ}&\text
{ for $t\in[0,1-(1/y)]$,}\\
W_t^{\circ}&\text{ for $t\in[1-(1/y),1].$}\endcases$$

We have
$$\rho\bigl(X_n\bigl(F_n(F_n^{-1}(\cdot))\bigr),
X_n\bigl(F_n(\tilde{F}_n^{-1}(\cdot))\bigr)\bigr)$$
 $$\qquad\qquad\qquad=\sup_{x\in[0,(1-\sqrt y)^2]}\vert
X_n\bigl(F_n(x)\bigr)-X_n\bigl(F_n((1-\sqrt y)^2)\bigr)\vert$$
$$\qquad\qquad\darrow
\sup_{x\in[0,(1-\sqrt y)^2]}\vert W_{F_y(x)}^{\circ}-W_{F_y((1-\sqrt
y)^2)}^{\circ}\vert=0$$ 
which implies
$$\rho\bigl(X_n\bigl(F_n(F_n^{-1}(\cdot))\bigr),
X_n\bigl(F_n(\tilde{F}_n^{-1}(\cdot))\bigr)\bigr)\parrow 0$$
Therefore (by Theorem 4.1 of [2])
$$X_n\bigl(F_n(F_n^{-1}(\cdot))\bigr)\darrow W_{F_y(F_y^{-1}(\cdot))}^{\circ}.$$

For $t<F_n(0)+\frac 1n$
$$X_n(t)=\sqrt{\frac
n2}\left(a_1^2\,\,\frac{\dsize\sum_{i=1}^{[nt]}z_i^2}{\dsize\sum_{\ell=1}^
{nF_n(0)}z_\ell^2}\quad-\quad \frac{[nt]}{n}\right)$$ $$=\quad
\frac{a_1^2}{\sqrt{F_n(0)}}\sqrt{\frac{nF_n(0)}{2}}\left
(\frac{\dsize\sum_{i=1}^{[nt]}z_i^2}{\dsize\sum_{\ell=1}^
{nF_n(0)}z_\ell^2}\quad-\quad\frac{[nt]}{nF_n(0)}\right)
\quad +\quad\tfrac{[nt]}{nF_n(0)}\sqrt{\tfrac
n2}(a_1^2-F_n(0))$$

Notice that $\sqrt{\frac n2}(a_1^2-F_n(0))=X_n(F_n(0))$.

For $t\in [0,1]$ let $\varphi_n(t)=\min(\frac{t}{F_n(0)},1)$,\
$\varphi(t)=\min(\frac{t}{1-(1/y)},1)$, and
$$Y_n(t)\quad=\quad\sqrt{\frac
n2}\left(\frac{\dsize\sum_{i=1}^{[nt]}z_i^2}{\dsize\sum_{\ell=1}^
{n}z_\ell^2}\quad- \quad\frac{[nt]}{n}\right).$$
Then $\varphi_n\parrow\varphi$ in $D_0\equiv\{x\in \underline D[0,1]:
x(1)\leq 1\}$ (see [2], p.144), and for $t<F_n(0)+\frac
1n$ $$Y_{nF_n(0)}(\varphi_n(t))\quad
=\quad\sqrt{\frac{nF_n(0)}{2}}\left(\frac{\dsize\sum_{i=1}^{[nt]}z_i^2}
{\dsize\sum_{\ell=1}^{nF_n(0)}z_\ell^2}\quad-\quad\frac{[nt]}{nF_n(0)}\right).$$

 For all $t\in [0,1]$ let
$$H_n(t)=\tfrac{a_1^2}{\sqrt{F_n(0)}}Y_{nF_n(0)}(\varphi_n(t))
+
X_n(F_n(0))\left(\tfrac{[nF_n(0)\varphi_n(t)]}{nF_n(0)}-1\right)+
X_n\bigl(F_n(F_n^{-1}(t))\bigr).$$ Then $H_n(t)=X_n(t)$ except on
intervals $[\frac mn,\frac{m+1}{n})$ where
$0<\lambda_m=\lambda_{m+1}$.  We will show $H_n\darrow W^{\circ}$ in
$D[0,1]$.

Let $\psi_n(t)=F_n(0)t$, $\psi(t)=(1-(1/y))t$,
 and 
$$V_n(t)=\frac1{\sqrt{2n}}\sum_{i=1}^{[nt]}(z_i^2-1).$$
Then $\psi_n\parrow \psi$ in
$D_0$ and
$$Y_n(t)=\frac{V_n(t)-\tfrac{[nt]}nV_n(1)}
{1+\sqrt{\tfrac2n}V_n(1)}.\leqno(2.6)$$

Since
$X_n\bigl(F_n(F_n^{-1}(\cdot))\bigr)$ and $V_n$ are independent we
have (using Theorems 4.4, 16.1 of [2])
$$\left(X_n\bigl(F_n(F_n^{-1}(\cdot))\bigr),V_n,\varphi_n,\psi_n
\right)\darrow\left(W_{F_y(F_y^{-1}(\cdot))}^{\circ},\overline
W,\varphi,\psi,\right)$$ where $\overline W$ is a Wiener process,
independent of $W^{\circ}$.  We immediately get ([2], p.145)
$$\left(X_n\bigl(F_n(F_n^{-1}(\cdot))\bigr),V_n\circ\psi_n,\varphi_n
\right)\darrow\left(W_{F_y(F_y^{-1}(\cdot))}^{\circ},\overline
W\circ\psi,\varphi\right).$$

Since $V_n(\psi_n(t))=\sqrt{F_n(0)}V_{nF_n(0)}(t)$, we have
$$\rho(V_n\circ\psi_n,\sqrt{1-(1/y)}V_{nF_n(0)})=\vert\sqrt{F_n(0)}-
\sqrt{1-(1/y)}\vert\sup_{t\in [0,1]}\vert V_{nF_n(0)}(t)\vert\parrow
0.$$ Therefore 
$$\left(X_n\bigl(F_n(F_n^{-1}(\cdot))\bigr),V_{nF_n(0)},\varphi_n
\right)\darrow\left(W_{F_y(F_y^{-1}(\cdot))}^{\circ},
\tfrac 1{\sqrt{1-(1/y)}}\overline
W\circ\psi,\varphi\right).\leqno(2.7)$$ 
Notice that	$\frac 1{\sqrt{1-(1/y)}}
\overline W\circ\psi$ is again a Wiener process, independent of $W^{\circ}$.

From (2.6) we have
$$Y_n(t)-(V_n(t)-tV_n(1))=V_n(t)\frac{t-\tfrac{[nt]}n+\sqrt{\tfrac2n}(t
V_n(1)-V_n(t))}{1+\sqrt{\tfrac2n}V_n(1)}.$$
Therefore $$\rho\bigl(Y_{nF_n(0)}(t),V_{nF_n(0)}(t)-tV_{nF_n(0)}(1)\bigr)
\parrow0.\leqno(2.8)$$

From (2.7), (2.8), and the fact that $W_t-tW_1$ is Brownian bridge 
it follows that
$$\left(X_n\bigl(F_n(F_n^{-1}(\cdot))\bigr),Y_{nF_n(0)},\varphi_n
\right)\darrow\left(W_{F_y(F_y^{-1}(\cdot))}^{\circ},\widehat
W^{\circ},\varphi\right)$$ where $\widehat W^{\circ}$ is another Brownian
bridge, independent of $W^{\circ}$.

The mapping $h:D[0,1]\times D[0,1]\times D_0\longrightarrow D[0,1]$
defined by $$h(x_1,x_2,z)=\sqrt{1-(1/y)}x_2\circ z+x_1(0)(z-1)+x_1$$ is
measurable, and is continuous on ${C[0,1]\times C[0,1]\times D\cap
C[0,1]}$.  Also, from (2.4) we have $a_1^2\parrow 1-(1/y)$.  Finally,
it is easy to verify $$\frac{[nF_n(0)\varphi_n(t)]}{nF_n(0)}\parrow
0$$
Therefore, we can conclude (using Theorem 4.1 and Corollary 1 of
Theorem 5.1 of [2]) 
$$H_n\darrow \sqrt{1-(1/y)}\widehat
W^{\circ}\circ\varphi+W_{1-(1/y)}^{\circ}(\varphi-1)+W_{F_y(F_y^{-1}(\cdot))}^{\circ}\equiv
H.$$ 

It is immediately clear that $H$ is a mean 0 Gaussian
process lying in $C[0,1]$. It is a routine matter to verify for
$0\leq s\leq t\leq1$ $$\text{\pe E}(H_sH_t)=s(1-t).$$

Therefore, $H$ is Brownian bridge.

We see that $\rho(X_n,H_n)$ is the same as the right hand side of
(2.3) except $i=1$ is excluded.  The arguments leading to (2.4) and
(2.5) ($2\leq i\leq t$) are exactly the same as before.  The fact that 
${\dsize\max_{2\leq i\leq t}}{m_i}/n\parrow 0$ follows from the case $y\leq 1$ since the
non-zero eigenvalues (including multiplicities) of $AA^T$ and $A^TA$
are identical for any rectangular $A$.  Thus
$$\rho(X_n,H_n)\parrow 0$$ and we have $X_n$ converging weakly to
Brownian bridge. $\square$ 
\vskip.1in

 \def\sqr#1#2{{\vcenter{\vbox{\hrule height.#2pt
\hbox{\vrule width.#2pt height#1pt \kern#1pt
\vrule width.#2pt}\hrule height.#2pt}}}}
\def\square{\mathchoice\sqr54\sqr54\sqr{2.1}3\sqr{1.5}3}

{\bf3. A new condition for weak convergence.}  In this section we establish two 
results on random elements of $D[0,b]$ needed for the proof of the 
Theorem 1.1.  In the following, we denote the modulus of continuity of 
$x\in D[0,b]$ by $w(x,\cdot)$:
$$w(x,\delta)=\sup_{|s-t|<\delta}|x(s)-x(t)|,\quad \delta\in(0,b].$$
To simplify the analysis we assume, for now, $b=1$.
\medskip
{\smc Theorem 3.1}.  {\sl Let $\{X_n\}$ be a sequence of random
elements of $D[0,1]$ whose probability measures satisfy the
assumptions of Theorem 15.5 of $[2]$, that is, $\{X_n(0)\}$ is tight,
and for every positive $\epsilon$ and $\eta$, there exists a
$\delta\in (0,1)$ and an integer $n_0$, such that, for all $n>n_0$,
{\pe P}$(w(X_n,\delta)\geq\epsilon)\leq\eta$.  If there exists a random
element $X$ with {\pe P}$(X\in C[0,1])=1$ and such that }
$$\left\{\int_0^1 t^rX_n(t)dt\right\}_{r=0}^\infty\darrow  
\left\{\int_0^1 t^rX(t)dt\right\}_{r=0}^\infty\quad\text{ as }
n\rightarrow\infty\leqno(3.1)$$ {\sl ($(\cal D)$ in (3.1) denoting weak
convergence on $R^\infty$), then $X_n\darrow X$.}
\smallskip
{\smc Proof}.  Note that the mappings 
$$x\longrightarrow \int_0^1 t^rx(t)dt$$ are continuous in $D[0,1]$.  
Therefore, by Theorems 5.1 and 15.5, $X_n\darrow X$ will follow if we
can show the distribution of $X$ is uniquely determined by the
distribution of
$$\left\{\int_0^1 t^rX(t)dt\right\}_{r=0}^\infty.\leqno(3.2)$$

Since the finite dimensional distributions of $X$ uniquely determine
the distribution of $X$, it suffices to show for any integer $m$ and numbers
$a_i,t_i$, ${i=0,1,\ldots,m}$ with $0=t_0<t_1<\cdots<t_m=1$, the 
distribution of 
$$\sum_{i=0}^m a_iX(t_i)\leqno(3.3)$$ is uniquely determined by the
distribution of (3.2).

Let $\{f_n\},f$ be uniformly bounded measurable functions on [0,1]
such that $f_n\rightarrow f$ pointwise as $n\rightarrow\infty$. 
Using the dominated convergence theorem we have
$$\int_0^1 f_n(t)X(t)dt\rightarrow\int_0^1 f(t)X(t)dt\quad\text{ as }n
\rightarrow\infty.\leqno(3.4)$$

Let $\epsilon>0$ be any number less than half the minimum distance between
the $t_i$'s.  Notice for the indicator function $I_{[a,b]}$ we have
the sequence of continuous ``ramp'' functions $\{R_n(t)\}$ with 
$$R_n(t)=\cases 1&t\in
[a,b],\\0&t\in[a-\frac1n,b+\frac1n]^c,\endcases$$
and linear on each of the sets {$[a-\frac1n,a]$},
{$[b,b+\frac1n]$}, satisfying $R_n\downarrow I_{[a,b]}$ as
$n\rightarrow\infty$.  Notice also that we can approximate any ramp
function uniformly on [0,1] by polynomials.  Therefore, using (3.4) for
polynomials appropriately chosen, we find that the distribution of 
$$\sum_{i=0}^{m-1}a_i\int_{t_i}^{t_i+\epsilon}X(t)dt+a_m\int_
{1-\epsilon}^1X(t)dt\leqno(3.5)$$ is uniquely determined by the distribution of
(3.2).  

Dividing (3.5) by $\epsilon$ and letting $\epsilon\rightarrow0$ we get
a.s. convergence to (3.3), (since $X\in C[0,1]$ with probability one), 
and we are done. $\square$
\medskip

{\smc Theorem 3.2.} {\sl Let $X$ be a random element of $D[0,1]$. 
Suppose there exists constants $B>0$, $\gamma\geq0$, $\alpha>1$, and a random
nondecreasing,
right-continuous function $F:[0,1]\rightarrow[0,B]$ such that, for
all $0\leq t_1\leq t_2\leq1$ and $\lambda>0$}
$$\text{\pe P}(|X(t_2)-X(t_1)|\geq\lambda)\leq\frac1{\lambda^{\gamma}}
\text{\pe E}[(F(t_2)-F(t_1))^{\alpha}].\leqno(3.6)$$
{\sl Then for every $\epsilon>0$ and $\delta$, an inverse of a positive
integer, we have}
$$\text{\pe P}(w(X,\delta)\geq3\epsilon)\leq\frac
{KB}{\epsilon^{\gamma}}\text{\pe E}\left[\max_{j<\delta^{-1}}
\left(F((j+1)\delta)-F(j\delta)\right)^{\alpha-1}\right],\leqno(3.7)$$
{\sl where $j$ ranges on positive integers, and $K$ depends only on $\gamma$ 
and $\alpha$.}

This theorem is proven by modifying the proofs of the first three theorems in
section 12 of [2].  It is essentially an extension of part of a result
contained in Theorem 12.3 of [2].  The original arguments, for the
most part, remain unchanged.  We will indicate only the specific
changes and refer the reader to [2] for details.  The extensions of
two of the theorems in [2] will be given below as lemmas.  However,
some definitions must first be given.

  Let $\xi_1,\ldots,\xi_m$ be random variables, and 
$S_k=\xi_1+\cdots+\xi_k$ ($S_0=0$).  Let $$M_m=\max_{0\leq k\leq
m}|S_k|\quad\text{ and}$$
$$M'_m=\max_{0\leq k\leq m}\min(|S_k|,|S_m-S_k|).$$ 
\smallskip

{\smc Lemma 3.1} (extends Theorem 12.1 of [2]). {\sl Suppose $u_1,
\ldots,u_m$ are non-negative random variables such that}
$$\text{\pe P}
(|S_j-S_i|\geq\lambda,|S_k-S_j|\geq\lambda)\leq\frac1{\lambda^
{2\gamma}}\text{\pe E}\biggl[\biggl(\sum_{i<\ell\leq k}u_\ell
\biggr)^{2\alpha}\biggr]<\infty,\quad\quad 0\leq i\leq j\leq k\leq m$$	 
{\sl for some $\alpha>\frac12$, $\gamma\geq0$, and for all
$\lambda>0$.  Then, for all $\lambda>0$ }
$$\text{\pe P}(M'_m\geq\lambda)\leq \frac
K{\lambda^{2\gamma}}\text{\pe E}[(u_1+\cdots+u_m)^{2\alpha}],\leqno(3.8)$$ where 
$K=K_{\gamma,\alpha}$ depends only on $\gamma$ and $\alpha$.
\smallskip

{\smc Proof}. We follow [2], p. 91. The constant $K$ is chosen
the same way and the proof proceeds by induction on $m$.  The
arguments for $m=1$ and 2 are the same, except, for the latter,
 $(u_1+u_2)^{2\alpha}$ is replaced by {\pe E}$(u_1+u_2)^{2\alpha}$.   
Assuming (3.8) is true for all integers less
than $m$, we find an integer $h$, $1\leq h\leq m$ such that
$$\frac{\text{\pe E}[(u_1+\cdots+u_{h-1})^{2\alpha}]}
{\text{\pe E}[(u_1+\cdots+u_m)^{2\alpha}]}\leq\frac12\leq
\frac{\text{\pe E}[(u_1+\cdots+u_h)^{2\alpha}]}
{\text{\pe E}[(u_1+\cdots+u_m)^{2\alpha}]},$$
the sum on the left hand side being 0 if $h=1$.

Since $2\alpha>1$, we have for all nonnegative $x$ and $y$
$$x^{2\alpha}+y^{2\alpha}\leq(x+y)^{2\alpha}.$$ We have then
$$\displaylines{\text{\pe E}[(u_{h+1}+\cdots+u_m)^{2\alpha}]\leq 
\text{\pe E}[(u_1+\cdots+u_m)^{2\alpha}]-
\text{\pe E}[(u_1+\cdots+u_h)^{2\alpha}]\hfill\cr
\hfill\leq\text{\pe E}[(u_1+\cdots+u_m)^{2\alpha}](1-\tfrac12)=
\tfrac12\text{\pe E}[(u_1+\cdots+u_m)^{2\alpha}].\cr}$$
Therefore, defining $U_1,U_2,D_1,D_2$ as in [2], we get the same
inequalities as in (12.30)-(12.33) ([2], p. 92) with $u^{2\alpha}$
replaced by {\pe E}$[(u_1+\cdots+u_m)^{2\alpha}]$.  The rest of the
proof follows exactly. $\square$
\medskip

{\smc Lemma 3.2} (extends Theorem 12.2 of [2]). {\sl If, for random
nonnegative $u_\ell$, there exists $\alpha>1$ and $\gamma\geq0$ such
that, for all $\lambda>0$}
$$\text{\pe P}(|S_j-S_i|\geq\lambda)\leq\frac1{\lambda^\gamma}
\text{\pe E}\biggl[\biggl(\sum_{i<\ell\leq j}u_\ell
\biggr)^{2\alpha}\biggr]<\infty,\quad 0\leq i\leq j\leq	m$$
{\sl then }
$$\text{\pe P}(M_n\geq\lambda)\leq
\frac{K'_{\gamma,\alpha}}{\lambda^\gamma}\text{\pe E}
[(u_1+\cdots+u_m)^{2\alpha}],\quad\quad\quad K'_{\gamma,\alpha}=
2^\gamma(1+K_{\frac12\gamma,\frac12\alpha}).$$
\smallskip

{\smc Proof}.  Following [2] we have for $0\leq i\leq j\leq
k\leq m$
$$\text{\pe P}
(|S_j-S_i|\geq\lambda,|S_k-S_j|\geq\lambda)\leq\text{\pe P}^{\frac12}
(|S_j-S_i|\geq\lambda)\text{\pe P}^{\frac12}
(|S_k-S_j|\geq\lambda)\leq\frac1{\lambda\gamma}\text{\pe E}
\biggl[\biggl(\sum_{i<\ell\leq k}u_\ell\biggr)^{2\alpha}\biggr],$$
so Lemma 3.1 is satisfied with constants $\frac12\gamma,\frac12\alpha$.  
The rest follows exactly as in [2], p. 94, with $(u_1+\cdots+u_m)^\alpha$ in
(12.46), (12.47) replaced by the expected value of the same quantity.
$\square$
\smallskip
We can now proceed with the proof of Theorem 3.2.  Following the
proof of Theorem 12.3 of [2] we fix positive integers $j<\delta^{-1}$ 
and $m$ and
define $$\xi_i=X\left(j\delta+\frac
im\delta\right)-X\left(j\delta+\frac{i-1}m\delta\right),\quad\quad
i=1,2,\ldots,m.$$  The partial sums of the $\xi_i$'s satisfy Lemma 3.2
with $$u_i=F\left(j\delta+\frac
im\delta\right)-F\left(j\delta+\frac{i-1}m\delta\right).$$
Therefore
$$\text{\pe P}\left({\max_{1\leq i\leq m}\left|{X\left(j\delta+\frac
im\delta\right)-X\left(j\delta\right)}\right|\geq\epsilon}\right)\leq
\frac K{\epsilon^\gamma}\text{\pe E}
[(F((j+1)\delta)-F(j\delta))^\alpha]\quad\quad
K=K'_{\gamma,\alpha}.$$

Since $X\in D[0,1]$ we have
$$\displaylines{\text{\pe P}\left(\sup_{j\delta\leq s\leq
(j+1)\delta}|X(s)-X(j\delta)|>\epsilon\right)\hfill\cr\hfill=\text{\pe P}
\left({\max_{1\leq i\leq m}\left|{X\left(j\delta+\frac
im\delta\right)-X\left(j\delta\right)}\right|>\epsilon}\text{ for all $m$
sufficiently large }\right)\cr\hfill\leq \lim_m\inf\text{\pe P}
\left({\max_{1\leq i\leq m}\left|{X\left(j\delta+\frac
im\delta\right)-X\left(j\delta\right)}\right|\geq\epsilon}\right)
\leq\frac K{\epsilon^\gamma}\text{\pe E}
[(F((j+1)\delta)-F(j\delta))^\alpha].\cr}$$
By considering a sequence of numbers approaching $\epsilon$ from
below we get from the continuity theorem
$$\text{\pe P}\left(\sup_{j\delta\leq s\leq
(j+1)\delta}|X(s)-X(j\delta)|\geq\epsilon\right)\leq 
\frac K{\epsilon^\gamma}\text{\pe E}
[(F((j+1)\delta)-F(j\delta))^\alpha].\leqno(3.9)$$
Summing both sides of (3.9) over all $j<\delta^{-1}$ and using the
corollary to Theorem 8.3 of [2] we get
$$\displaylines{\text{\pe P}(w(X,\delta)\geq3\epsilon)\leq
\frac K{\epsilon^\gamma}\text{\pe E}\biggl[\sum_{j<\delta^{-1}}
(F((j+1)\delta)-F(j\delta))^\alpha\biggr]\hfill\cr\hfill
\leq\frac K{\epsilon^\gamma}\text{\pe E}\biggl[
\max_{j<\delta^{-1}}(F((j+1)\delta)-F(j\delta))^{\alpha-1}(F(1)-F(0))
\biggr]\cr\hfill\leq\frac{KB}{\epsilon^\gamma}\text{\pe E} 
\biggl[\max_{j<\delta^{-1}}(F((j+1)\delta)-F(j\delta))^{\alpha-1}\biggr],
\cr}$$ and we are done. $\square$

For general $D[0,b]$ we simply replace (3.1) by  
$$\left\{\int_0^b t^rX_n(t)dt\right\}_{r=0}^\infty\darrow  
\left\{\int_0^b t^rX(t)dt\right\}_{r=0}^\infty\quad\text{ as }
n\rightarrow\infty\leqno(3.10)$$ and (3.7) by
$$\text{\pe P}(w(X,b\delta)\geq3\epsilon)\leq\frac
{KB}{\epsilon^{\gamma}}\text{\pe E}\left[\max_{j<\delta^{-1}}
\left(F(b(j+1)\delta)-F(bj\delta)\right)^{\alpha-1}\right],\leqno(3.11)$$
$j$ and $\delta^{-1}$ still positive integers.
 \vskip.1in
 
{\bf4. Completing the proof.} We finish up by verifying the conditions
of Theorem 3.1.
\medskip
{\smc Theorem 4.1}. {\sl Let $\text{\pe E}(v_{11})=0$, $\text{\pe E}
(v_{11}^2) = 1$, and $\text{\pe E}(v_{11}^4)<\infty$. Suppose the sequence of 
vectors $\{\vec x_n\}$, $\vec x_n=(x_{n1},x_{n2},\ldots,x_{nn})^T$, 
$\|\vec x_n\| = 1$ satisfies}
$$\sum_{i=1}^n x_{ni}^4\,\longrightarrow\,0\quad\text{as
$n\rightarrow\infty$}.\leqno(4.1)$$ {\sl Then (1.5) holds.}

\smallskip
{\smc Proof}. Let $\bar v_{ij}=\bar v_{ij}(n)=v_{ij}I_{(|v_{ij}|
\leq n^{1/4})}-
\text{\pe E}(v_{ij}I_{(|v_{ij}|\leq n^{1/4})})$, and let $\bar M_n=\frac 1s\bar
V_n\bar V_n^T$, where $\bar V_n=(\bar v_{ij})$.  We have $\text{\pe E}(\bar
v_{11})=0$, $\text{\pe E}(\bar v_{11}^2)\rightarrow 1$, and $\text{\pe E}(\bar
v_{11}^4)\rightarrow \text{\pe E}(v_{11}^4)$ as $n\rightarrow\infty$.

The main part of the proof in [11] establishing (1.5) (under the
additional assumption $\text{\pe E}(v_{11}^4)=3$) relies on a 
multidimensional version of the method of moments, together with the fact
that (1.5) holds in the Wishart case.  It is shown that for 
any integer $m\geq 2$ and positive integers $r_1,r_2,\ldots,r_m$, the 
asymptotic behavior of 
$$\displaylines{(4.2)\quad n^{m/2}\text{\pe E}\left[{\left(\vec x_n^T\bar M_n^{r_1}\vec x_n-
\text{\pe E}(\vec x_n^T\bar M_n^{r_1}\vec x_n)\right)}{\left(\vec x_n^T\bar
M_n^{r_2}\vec x_n-
\text{\pe E}(\vec x_n^T\bar M_n^{r_2}\vec x_n)\right)}\right.\hfill\cr
\hfill\left.\cdots{\left(\vec x_n^T\bar
M_n^{r_m}\vec x_n-
\text{\pe E}(\vec x_n^T\bar M_n^{r_m}\vec x_n)\right)}\right]\cr}$$
depends only on $\text{\pe E}(\bar v_{11}^2)$ and $\text{\pe E}
(\bar v_{11}^4)$.  Using the fact that (4.2) converges to the appropriate 
limit when $v_{11}$ is N(0,1), (1.5) follows when $\text{\pe E}(v_{11}^4)=3$.

This is the only place in the proof in [11] that refers to the value of
$\text{\pe E}(v_{11}^4)$, the remaining arguments depending on this value 
only to the extent of it being finite, and thus apply to the present case.  
Therefore, we will be done if we can show that (4.1) implies the asymptotic 
behavior of (4.2) depends only on $\text{\pe E}(\bar v_{11}^2)$. Although it will be 
necessary to repeat some of the discussion in the original proof, we refer 
the reader to [11] for specific details.

We have (dropping the dependency of $n$ on the components of $\vec
x_n$)
$$\left(\frac{s^{r_1+\cdots+r_m}}{n^{m/2}}\right)\times(4.2)\,\,=$$
$$\displaylines{\text{(4.3) }\sum
x_{i^1}x_{j^1}\cdots x_{i^m}x_{j^m}\text{\pe E}\left[(\bar
v_{i^1k_1^1}\bar v_{i_2^1k_1^1}\cdots\bar v_{i_{r_1}^1k_{r_1}^1}
\bar v_{j^1k_{r_1}^1}-\text{\pe E}(\bar
v_{i^1k_1^1}\bar v_{i_2^1k_1^1}\cdots\bar v_{i_{r_1}^1k_{r_1}^1}
\bar v_{j^1k_{r_1}^1}))\right.\hfill\cr
\hfill\cdots\left.(\bar
v_{i^mk_1^m}\bar v_{i_2^mk_1^m}\cdots\bar v_{i_{r_m}^mk_{r_m}^m}
\bar v_{j^mk_{r_m}^m}-\text{\pe E}(\bar
v_{i^mk_1^m}\bar v_{i_2^mk_1^m}\cdots\bar v_{i_{r_m}^mk_{r_m}^m}
\bar v_{j^mk_{r_m}^m}))\right],\cr}$$
where the sum is over $i^1,j^1,i_2^1,\ldots,i_{r_1}^1,k_1^1,\ldots,k_{r_1}^1
,\ldots,i^m,j^m,i_2^m,\ldots,i_{r_m}^m,k_1^m,\ldots,k_{r_m}^m$.

We consider one of the ways the two set of indices 
\flushpar$I\equiv\{i^1,j^1,i_2^1,\ldots,i_{r_1}^1,\ldots,
i^m,j^m,i_2^m,\ldots,i_{r_m}^m\}$, and
$K\equiv\{k_1^1,\ldots,k_{r_1}^1,\ldots,k_1^m,\ldots,k_{r_m}^m\}$
can each be partitioned.  Associated with the two partitions are the
terms in (4.3) (for $n$ large) where indicies are equal in value if and
only if they belong to the same class.  We consider only those
partitions corresponding to terms in (4.3) that contribute a 
non-negligible amount to (4.2) in the limit.
Let $\ell$ be the number of classes of $I$ indices containing only
one element from $J\equiv\{i^1,j^1,i^2,j^2,\ldots,i^m,j^m\}$.  Let
$d$ denote the number of classes of $I$ indices containing no
elements from $J$, plus the number of classes of $K$ indices.  Then
the contribution to (4.3) of those terms associated with the two
partitions is bounded in absolute value by
$$Cn^{(\ell/2)+d}\text{\pe E}\left(|\bar v_{i^1k_1^1},\ldots,\bar v_{j^1k_{r_1}^1},
\ldots,\bar v_{i^mk_1^m},\ldots,\bar
v_{j^mk_{r_m}^m}|\right)\leqno(4.4)$$
the expected value being one of those associated with the two
partitions.  We can write this expected value in the form  
$$A_{a_1b_1}^1A_{a_2b_2}^2\cdots A_{a_{r'}b_{r'}}^{r'}$$
where $A_{a_jb_j}^j$ corresponds to $\bar v_{a_jb_j}$ appearing in
(4.4), so that if $\bar v_{a_jb_j}$ appears $t$ times, then 
$A_{a_jb_j}^j=\text{\pe E}(|\bar v_{a_jb_j}^t|)$.  There are $r'$ distinct
elements of $\bar V_n$ involved in (4.4).  For each ordered pair
$(a_j,b_j)$ either $a_j$ or $b_j$ will be repeated in at least one
other ordered pair (see [11]).
We say that $a_j$ or $b_j$ is \underbar{free} if it does  
not appear in any other ordered pair.

From [11] it was argued that $d=r_1+\cdots+r_m-(m/2)-(\ell/2)$, 
for each $j$ $A_{a_jb_j}^j=\text{\pe E}(\bar v_{11}^2)$ or $\text{\pe E}(\bar v_{11}^4)$, and, for our
purposes, we can assume without loss of generality that each $I$
class containing no element from $J$ has at least two elements, and
any $A_{a_jb_j}^j$ for which $b_j$ is free involves 
$\bar v_{p^tq^t}$'s for at least two different $t$'s.  We continue
the argument from this point.

We immediately conclude that the number of $I$ classes containing 
no elements from $J$, and the number of $K$ classes are,
respectively, bounded by $\frac12(r_1+\cdots+r_m-m)$ and 
$\frac12(r_1+\cdots+r_m)$.  We can also conclude that $\ell=0$, since
any $I$ class containing only one element from $J$ implies an element
of $\bar V_n$ apearing in (4.4) an odd number of times, which is
impossible.  Therefore, we can conclude that there are
$\frac12(r_1+\cdots+r_m)$ $K$ classes, and
$\frac12(r_1+\cdots+r_m-m)$ $I$ classes containing no
elements from $J$, which further implies no $I$ class 
exists containing an element from $J$ and an element 
from $(I-J)$.  Therefore, the $I$ classes split up into separate
$J$ and $I-J$ classes, and
each $K$ class and $I-J$ class
consists of two elements.

Now, if $A_{a_jb_j}^j=\text{\pe E}(\bar v_{11}^4)$	for some $j$, then either
$a_j$ or $b_j$ must be associated with a class containing at least
three elements, forcing $a_j$ to be associated with a $J$ class
consisting of at least four elements.  Therefore, the sum of the terms in
(4.3) associated with the two partitions, divided by 
$(s^{r_1+\cdots+r_m})/n^{m/2}$, will be bounded by 
$$C\times\sum_{i=1}^nx_i^4$$   
and because of (4.1) these terms will not contribute anything to (4.2) 
in the limit.  We conclude that the asymptotic behavior of (4.2) 
depends only on $\text{\pe E}(\bar v_{11}^2)$, and we are done. $\square$

Let $R_+$ denote the nonnegative reals and $\cal B_+$, 
$\cal B_+^4$ denote the Borel $\sigma$-fields on, respectively, $R_+$ and 
$R_+^4$.  For any $n\times n$ symmetric, nonnegative definite matrix $B$
and any $A\in {\cal B}_+$, let 
$P^B(A)$ denote the projection matrix on the
subspace of $R^n$ spanned by the eigenvectors of $B$ having
eigenvalues in $A$ (the collection of projections 
$\{P^B((-\infty,a]):\,a\in R\}$ is usually referred to as the
spectral family of $B$).  We have tr$P^B(A)$ equal to the number of 
eigenvalues of $B$ contained in $A$.  If $B$ is random, then it is
straightforward to verify the following facts:
\medskip
\flushpar a) For every $\vec x_n\in R^n$, $\|\vec x_n\|=1$, 
$\vec x_n^TP^B(\cdot)\vec x_n$ is a
random probability measure on $R_+$ placing mass on the eigenvalues
of $B$.

\medskip
 \flushpar b) For any four elements $P_{i_1j_1}^B(\cdot)$,
$P_{i_2j_2}^B(\cdot)$, $P_{i_3j_3}^B(\cdot)$, $P_{i_4j_4}^B(\cdot)$ of
$P^B(\cdot)$, the function defined on rectangles $A_1\times A_2\times 
A_3\times A_4\in \cal
B_+^4$ by
$$\text{\pe E}\bigl(P_{i_1j_1}^B(A_1)P_{i_2j_2}^B(A_2)P_{i_3j_3}^B(A_3)
P_{i_4j_4}^B(A_4)\bigr)\leqno(4.5)$$ generates a signed measure
$m_n^B=m_n^{B,(i_1,j_1,\ldots,i_4,j_4)}$ on $(R_+^4,\cal B_+^4)$
such that $|m_n^B(A)|\leq 1$ for every $A\in \cal B_+^4$.
\smallskip
When $B=M_n$ we also have
\medskip 
\flushpar c) For any $A\in \cal B_+$ the distribution of $P^{M_n}(A)$
is invariant under permutation transformations, that is, 
$P^{M_n}(A)\sim OP^{M_n}(A)O^T$ for any permutation matrix $O$ (use
the fact that $P^B(\cdot)$ is uniquely determined by
$\{B^r\}_{r=1}^{\infty}$ along with $OP^B(\cdot)O^T=P^{OBO^T}(\cdot)$
and $\{M_n^r\}_{r=1}^{\infty}\sim \{(OM_nO^T)^r\}_{r=1}^{\infty}$).
\medskip
\flushpar d) For $0\leq x_1\leq x_2$ 
$$\tfrac1nP^{M_n}([0,x_1])=F_n(x_1),$$
$$X_n(F_n(x_1))=\sqrt{\tfrac n2}\bigl(\vec x_n^TP^{M_n}([0,x_1])\vec
x_n-\tfrac1n\text{tr}(P^{M_n}([0,x_1]))\bigr),\quad\text{and}$$ 
$$X_n(F_n(x_2))-X_n(F_n(x_1))=\sqrt{\tfrac n2}\bigl(\vec
x_n^TP^{M_n}((x_1,x_2])
\vec x_n-\tfrac1n\text{tr}(P^{M_n}((x_1,x_2]))\bigr).$$
\medskip

{\smc Lemma 4.1}. {\sl Assume $v_{11}$ is symmetric.
If one of the indices $i_1,j_1,\ldots,i_4,j_4$ appears
an odd number of times, then $m_n^{M_n}\equiv 0$.}
\medskip

{\smc Proof}. Assume first that $v_{11}$ is 
bounded.  Then $\lambda_{\max}$ is bounded which implies $m_n^{M_n}$
has bounded support.  Therefore, $m_n^{M_n}$ is uniquely determined
by its mixed moments ([7], pp. 97-102).  It is straightforward to show  
these moments can be expressed as
$$\text{\pe E}\bigl(({M_n}^{r_1})_{i_1j_1}({M_n}^{r_2})_{i_2j_2}({M_n}^{r_3})_{i_3j_3}
({M_n}^{r_4})_{i_4j_4}\bigr)\leqno(4.6)$$
for arbitrary nonnegative integers $r_1,\ldots,r_4$.  If, say
$r_1=0$ and $i_1\neq j_1$, then obviously (4.6) is zero.  We can assume
then a positive power $r_{\ell}$ for which $i_{\ell}\neq j_{\ell}$.  Upon
expanding (4.6) as a sum of expected values of products of entries 
of $V_n$ we find each product contains a $v_{ij}$ appearing an odd
number of times.  Therefore, all mixed moments of $m_n^{M_n}$ are zero, 
implying $m_n^{M_n}\equiv 0$.

For arbitrary symmetric $v_{11}$ we truncate.  For any $c>0$ let 
$v_{ij}^c=v_{ij}I_{(|v_{ij}|\leq c)}$, and
$M_n^c=\frac1sV_n^c{V_n^c}^T$, where $V_n^c=(v_{ij}^c)$.  For any 
realization of the $v_{ij}$'s and $A_1\times A_2\times A_3\times 
A_4\in \cal B_+^4$
$$P_{i_1j_1}^{M_n^c}(A_1)P_{i_2j_2}^{M_n^c}(A_2)
P_{i_3j_3}^{M_n^c}(A_3)P_{i_4j_4}^{M_n^c}(A_4)
\quad=\quad P_{i_1j_1}^{M_n}(A_1)P_{i_2j_2}^{M_n}(A_2)
P_{i_3j_3}^{M_n}(A_3)P_{i_4j_4}^{M_n}(A_4)$$
for $c$ large enough.  Therefore, by the dominated convergence
theorem we have (4.5) (with $B=M_n$) equal to zero for all rectangles in 
$\cal B_+^4$, implying $m_n^{M_n}\equiv 0$. $\square$
\medskip

{\smc Theorem 4.2}. {\sl Assume $v_{11}$ is symmetric
and $\vec x_n=(\pm\frac1{\sqrt n},\pm\frac1{\sqrt
n},\ldots,\pm\frac1{\sqrt n})^T$.  Let $G_n(x)=4F_n(x)$.  Then }
$$\text{\pe E}\bigl((X_n(F_n(0)))^4\bigr)\leq 
\text{\pe E}\bigl((G_n(0))^2\bigr),\leqno(4.7)$$
{\sl and for any $0\leq x_1\leq x_2$}
$$\text{\pe E}\bigl((X_n(F_n(x_2))-X_n(F_n(x_1)))^4\bigr)\leq 
\text{\pe E}\bigl((G_n(x_2)-G_n(x_1))^2\bigr).\leqno(4.8)$$
\medskip

{\smc Proof.}  With $A=\{0\}$ in (4.7), $A=(x_1,x_2]$ in (4.8)
we use d) to find the left
hand sides of (4.7) and (4.8) equal to
$$\frac1{4n^2}\text{\pe E}\biggl(\sum_{i\neq
j}\gamma_{ij}P_{ij}^{M_n}(A)\biggr)^4\leqno(4.9)$$
where $\gamma_{ij}=\text{sgn}(\vec x_n)_i(\vec x_n)_j$.   For the
remainder of the argument we simplify the notation by supressing the
dependence of the projection matrix on $M_n$ and $A$.  Upon expanding
(4.9) we use c) to combine identically distributed factors, and 
Lemma 4.1 to arrive at
$$\displaylines{(4.10)\quad(4.9)\,=\,\tfrac{(n-1)}n\bigl(12(n-2)
\text{\pe E}(P_{12}^2P_{13}^2)+3(n-2)(n-3)\text{\pe E}(P_{12}^2P_{34}^2)\hfill\cr\hfill
+12(n-2)(n-3)\text{\pe E}(P_{12}P_{23}P_{34}P_{14})+2\text{\pe E}(P_{12}^4)\bigr).\cr}$$

We can write the second and third expected values in (4.10) in terms of
the first expected value and expected values involving $P_{11}$,
$P_{22}$, and $P_{12}$
by making further use of c) and the fact that P
is a projection matrix (i.e., $P^2=P$).  For example, we take the
expected value of both sides of the identity
$$P_{12}P_{23}\biggl(\sum_{j\geq
4}P_{3j}P_{1j}\,+P_{31}P_{11}+P_{32}P_{12}+P_{33}P_{13}\biggr)=P_{12}P_{
23}P_{13}$$ and get
$$(n-3)\text{\pe E}(P_{12}P_{23}P_{34}P_{14})+2\text{\pe E}(P_{11}P_{12}P_{23}P_{31})
+\text{\pe E}(P_{12}^2P_{13}^2)=\text{\pe E}(P_{12}P_{23}P_{13}).$$ Proceeding in the same
way we find $$(n-2)\text{\pe E}(P_{11}P_{12}P_{23}P_{31})+\text{\pe E}(P_{11}^2P_{12}^2)
+\text{\pe E}(P_{11}P_{22}P_{12}^2)=\text{\pe E}(P_{11}P_{12}^2)\quad\text{ and}$$
$$(n-2)\text{\pe E}(P_{12}P_{23}P_{13})+2\text{\pe E}(P_{11}P_{12}^2)=\text{\pe E}(P_{12}^2).$$
Therefore $$\displaylines{(n-2)(n-3)\text{\pe E}(P_{12}P_{23}P_{34}P_{14})=\text{\pe E}(P_{12}^2)+
2\text{\pe E}(P_{11}P_{22}P_{12}^2)\hfill\cr\hfill
+2\text{\pe E}(P_{11}^2P_{22}^2)-(n-2)\text{\pe E}(P_{12}^2P_{13}^2)
-4\text{\pe E}(P_{11}P_{12}^2).\cr}$$

Since $P_{11}\geq \max(P_{11}P_{22},P_{11}^2)$ and $P_{12}^2\leq 
P_{11}P_{22}$ (since $P$ is nonnegative definite) we have
$$(n-2)(n-3)\text{\pe E}(P_{12}P_{23}P_{34}P_{14})\leq \text{\pe E}(P_{11}P_{22})-
(n-2)\text{\pe E}(P_{12}^2P_{13}^2).$$

Similar arguments will yield 
$$(n-3)\text{\pe E}(P_{12}^2P_{34}^2)+2\text{\pe E}(P_{12}^2P_{13}^2)+\text{\pe E}(P_{12}^2P_{33}^2)=
\text{\pe E}(P_{12}^2P_{33})\quad\text{ and}$$
$$(n-2)\text{\pe E}(P_{12}^2P_{13}^2)+\text{\pe E}(_{12}^4)+\text{\pe E}(P_{11}^2P_{12}^2)=\text{\pe E}(P_{11}P_{12}^2).
$$  After multiplying the first equation by $n-2$ and adding it to the
second, we get
$$\displaylines{(n-2)(n-3)\text{\pe E}(P_{12}^2P_{34}^2)+3(n-2)\text{\pe E}(P_{12}^2P_{13}^2)
\hfill\cr\hfill=(n-2)\text{\pe E}(P_{12}^2P_{33})-(n-2)\text{\pe E}(P_{12}^2P_{33}^2)+\text{\pe E}(P_{11}P_{12}^2)-\text{\pe E}(P_
{11}^2P_{12}^2)-\text{\pe E}(P_{12}^4)\cr
\hfill=\text{\pe E}(P_{11}P_{22})+\text{\pe
E}(P_{11}^2P_{22}^2)-2\text{\pe
E}(P_{11}P_{22}^2)-\text{\pe E}(P_{12}^4)
\leq \text{\pe E}(P_{11}P_{22})-\text{\pe E}(P_{12}^4).\cr}$$

Combining the above expressions we obtain
$$\quad(4.9)\,\leq\,15\frac{(n-1)}n\text{\pe E}(P_{11}P_{22}).$$
Therefore, using c) and d), we get 
$$(4.9)\leq\tfrac{15}{n^2}\text{\pe E}({\tsize\sum_{i\neq
j}}P_{ii}P_{jj})\leq\text{\pe E}((4\tfrac1n\text{tr}P)^2) 
=\cases\text{\pe E}\bigl((G_n(0))^2\bigr)&\text{ for $A=\{0\}$},\\
\text{\pe E}\bigl((G_n(x_2)-G_n(x_1))^2\bigr)&\text{ for
$A=(x_1,x_2]$},\endcases $$
and we are done. $\square$
\medskip

We can now complete the proof of Theorem 1.1.  We may assume 
{\pe E}$(v_{11}^2)=1$.  Choose any $b>(1+\sqrt y)^2$.
We have (1.3) and, by Theorem 4.1, (1.5), which imply
$$\left\{\int_0^b x^rX_n(F_n(x))dx\right\}_{r=0}^\infty\darrow  
\left\{\int_0^b x^rW_x^ydx\right\}_{r=0}^\infty\quad\text{ as }
n\rightarrow\infty,$$  so that (3.10) is satisfied.  By Theorems 3.2 and 4.2  
we have, for any $n\geq16$, (3.11) with $X=X_n(F_n(\cdot))$, $F=4F_n$, 
 $B=4$, $\gamma=4$, and $\alpha=2$.  From (1.2) and Theorem 5.1 of
[2] we have for every $\delta\in(0,b]$
$$w(F_n,\delta)\parrow w(F_y,\delta)\quad\text{ as }n\rightarrow\infty.$$
Since $F_y$ is continuous on $[0,\infty)$, we apply the dominated
convergence theorem to the right hand side of (3.11) and find that,
for every $\epsilon>0$, {\pe
P}$(w(X_n(F_n(\cdot)),\delta)\geq\epsilon)$ can be made arbitrarily
small for all $n$ sufficiently large by choosing $\delta$
 appropriately.
Therefore, by Theorem 3.1, $X_n(F_n(\cdot))\darrow
W_{F_y(\cdot)}^\circ$ in
$D[0,b]$, which implies $X_n(F_n(\cdot))\darrow W_{F_y(\cdot)}^\circ$ in
$D[0,\infty)$, and by Theorem 2.1 we conclude that $X_n\darrow
W^\circ$. $\square$
\vskip.2in
{\bf Acknowledgment.} I would like to thank the Associate Editor and the 
referee for their
helpful suggestions.
\vfill\eject
\centerline{REFERENCES}
\medskip
\item{[2]} Billingsley, P. (1968). {\sl Convergence of Probability
Measures}. Wiley, New York.
\item{[3]}} Billingsley, P. (1995) {\sl Probability and Measure} Third Edition. Wiley,
New York.
\item{[1]}  Billingsley, P. (1999). {\sl Convergence of Probability
Measures} Second Edition. Wiley, New York.
\item{[4]} Grenander, U. and Silverstein, J.W. (1977). Spectral
analysis of networks with random topologies. {\sl SIAM J. Appl. Math.
\bf37} 499-519.
\item{[5]} Jonsson, D. (1982). Some limit theorems for the
eigenvalues of a sample covariance matrix. {\sl J. Multivariate Anal.
\bf12} 1-38.
\item{[6]} Lindvall, T. (1973). Weak convergence of probability
measures and random functions in the function space $D[0,\infty)$).
{\sl J. Appl. Probab. \bf10} 109-121.
\item{[8]} Silverstein, J.W. (1979). On the randomnes of eigenvectors
generated from networks with random topologies. {\sl SIAM J. Appl.
Math. \bf37} 235-245.
\item{[9]} Silverstein, J.W. (1981). Describing the behavior of
random matrices using sequences of measures on orthogonal groups.
 {\sl SIAM J. Math. Anal. \bf12} 274-281.
\item{[10]} Silverstein, J.W. (1984). Some limit theorems on the
eigenvectors of large dimensional sample covariance matrices. 
 {\sl J. Multivariate Anal. \bf15} 295-324.
\item{[11]} Silverstein, J.W. (1989). On the eigenvectors of large
dimensional sample covariance matrices. {\sl J. Multivariate Anal.}
{\bf 30} 1-16.
\item{[14]} Silverstein, J.W. (1990) Weak convergence of random functions defined by the
eigenvectors of sample covariance matrices.  {\sl Ann. Probab. \bf18 1174-1194
\item{[12]} Wachter, K.W. (1978). The strong limits of random matrix
spectra for sample matrices of independent elements. {\sl Ann.
Probab. \bf6} 1-18.
\item{[13]} Yin, Y.Q. (1986). Limiting spectral distribution for a
class of random matrices. {\sl J. Multivariate Anal. \bf20} 50-68.
\item{[14]} Yin, Y.Q., Bai, Z.D., and Krishnaiah, P.R. (1988). On limit of
the largest eigenvalue of the large dimensional sample covariance
matrix. {\sl Probab. Th. Rel. Fields \bf78} 509-521.
\vskip.5in
\hskip2.5in {\smc Department of Mathematics}\smallskip
\hskip2.5in {\smc Box 8205}\smallskip
\hskip2.5in{\smc North Carolina State University}\smallskip
\hskip2.5in{\smc Raleigh, North Carolina 27695-8205}
\bye